# Singular dual pairs

Juan–Pablo Ortega

Institut Nonlinéaire de Nice

Centre National de la Recherche Scientifique

1361, route des Lucioles

F-06560 Valbonne, France

Juan-Pablo.Ortega@inln.cnrs.fr

February 19, 2002

## Abstract

We generalize the notions of dual pair and polarity introduced by S. Lie [Lie90] and A. Weinstein [W83] in order to accommodate very relevant situations where the application of these ideas is desirable. The new notion of polarity is designed to deal with the loss of smoothness caused by the presence of singularities that are encountered in many problems in Poisson and symplectic geometry. We study in detail the relation between the newly introduced dual pairs, the quantum notion of Howe pair, and the symplectic leaf correspondence of Poisson manifolds in duality. The dual pairs arising in the context of symmetric Poisson manifolds are treated with special attention. We show that in this case and under very reasonable hypotheses we obtain a particularly well behaved kind of dual pairs that we call *von Neumann pairs*. Some of the ideas that we present in this paper shed some light on the *optimal momentum maps* introduced in [OR02a].

# Contents

# 1 Introduction

The notion of dual pair, introduced by A. Weinstein in [W83], is of central importance in the context of Poisson geometry. Let $(M, \omega)$ be a symplectic manifold, $(P_1, \{\cdot, \cdot\}_{P_1})$ and $(P_2, \{\cdot, \cdot\}_{P_2})$ be two Poisson manifolds, and $\pi_1 : M \to P_1$ and $\pi_2 : M \to P_2$ be two Poisson surjective submersions. The diagram is called a **dual pair** if the Poisson subalgebras $\pi_1^* C^\infty(P_1)$ and $\pi_2^* C^\infty(P_2)$ centralize each other with respect to the Poisson structure in $M$ associated to the symplectic form $\omega$. This notion has its origins in the study of group representations arising in quantum mechanics. In this direction, we have the works of Howe [Howe89], Kashiwara and Vergne [KV78], Sternberg and Wolf [SW78], and Jakobsen and Vergne [JV79], which justify why we will refer to the previously defined dual pairs as **Howe pairs**.

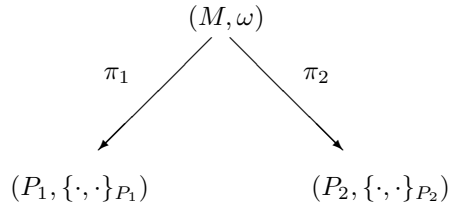

Already in 1890, S. Lie (see [Lie90] and §7 in [W83]) devised a method to construct Howe pairs using the notion of **polarity,** that we briefly describe: let $D$ be an integrable regular distribution on the symplectic manifold $(M, \omega)$ that is everywhere the span of locally Hamiltonian vector fields. Under these circumstances the space of leaves $M/D$ is a Poisson manifold and the canonical projection $\pi_D : M \to M/D$ is a Poisson surjective submersion [MR86]. Let now $D^\omega$ be the **polar** distribution to $D$, defined by

$$D^\omega(m) := \{v \in T_m M \mid \omega(m)(v, w) = 0 \quad \text{for all} \quad w \in D(m)\}.$$

A simple verification shows that $D^\omega$ is smooth and integrable. If we assume that the corresponding space of leaves $M/D^\omega$ is a regular quotient manifold and denote by $\pi_{D^\omega} : M \to M/D^\omega$ the canonical projection, then the diagram $M/D \xleftarrow{\pi_D} M \xrightarrow{\pi_{D^\omega}} M/D^\omega$ is a Howe pair. Moreover, $\ker T\pi_D$ and $\ker T\pi_{D^\omega}$ are symplectically orthogonal distributions. This remark motivates the following definition: the diagram $(P_1, \{\cdot, \cdot\}_{P_1}) \overset{\pi_1}{\leftarrow} (M, \omega) \overset{\pi_2}{\rightarrow} (P_2, \{\cdot, \cdot\}_{P_2})$ is called a ***Lie–Weinstein dual pair*** when $\ker T\pi_1$ and $\ker T\pi_2$ are symplectically orthogonal distributions. Every Lie–Weinstein dual pair where the submersions $\pi_1$ and $\pi_2$ have connected fibers is a Howe pair (see Corollary 5.4).

The geometries underlying two Poisson manifolds forming a Lie–Weinstein pair are very closely related. For instance, if the fibers of the submersions $\pi_1$ and $\pi_2$ are connected, the symplectic leaves of $P_1$ and $P_2$ are in bijection [KKS78, W83, Bl01] and, for any $m \in M$, the transverse Poisson structures on $P_1$ and $P_2$ at $\pi_1(m)$ and $\pi_2(m)$, respectively, are anti–isomorphic [W83].

Apart from the already mentioned studies on representation theory, dual pairs occur profusely in finite and infinite dimensional classical mechanics (see for instance [MRW84, MW83], and references therein). Another intimately related concept that we will not treat in our study is that of the ***Morita equivalence*** of two Poisson manifolds [Mor58, Xu91]. Nice presentations of the classical theory of dual pairs can be found in [CaWe99] and in [Land98].

In this paper we will pay special attention to the dual pairs that appear in the reduction of Poisson symmetric systems. We introduce this situation with a very simple example: let $(M, \omega)$ be a symplectic manifold and $G$ be a Lie group acting freely, canonically, and properly on $M$. Suppose that this action has a momentum map $\mathbf{J} : M \to \mathfrak{g}^*$ associated. If we denote by $\mathfrak{g}_{\mathbf{J}}^*$ the image of $M$ by $\mathbf{J}$, it is easy to check that the diagram is a Lie–Weinstein dual pair, and consequently, if the group is connected and $\mathbf{J}$ has connected fibers, a Howe pair. However, in most cases of practical interest, the freeness assumption on the $G$–action is not satisfied hence it is worth studying the impact of dropping this condition in the duality between $M/G$ and $\mathfrak{g}_{\mathbf{J}}^*$. When the action is not free $M/G$ and $\mathfrak{g}_{\mathbf{J}}^*$ still form a Lie–Weinstein dual pair in a generalized sense since, even though $M/G$ is not a smooth manifold anymore, the tangent space to the $G$–orbits (fibers of $\pi$) and $\ker T\mathbf{J}$ are symplectically orthogonal. The question now is: do they still form a Howe pair in the presence of the connectedness hypotheses? or, in other words: do the $G$–invariant functions $C^\infty(M)^G$ in $M$ and the ***collective functions*** $\mathbf{J}^*C^\infty(\mathfrak{g}^*)$ centralize each other? This question has deserved much attention due in part to physical motivations [GS80]. Guillemin and Sternberg conjectured in [GS84] that the answer to our question was affirmative for any compact Hamiltonian group action, and they proved it for toral actions. However, Lerman gave in [L89] a counterexample to this conjecture that showed the first indications of the great complexity underlying the relation between Lie–Weinstein and Howe pairs in the case of non free actions. This relation, that may eventually become very sophisticated, has been the subject of studies of great interest. See for instance [Kar93, KL97, Knop90, knop97], and references therein.

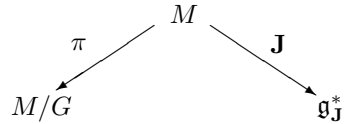

Another notion that breaks down in the absence of regularity hypotheses is that of polarity. In the previous paragraphs we mentioned that the polar distribution to an integrable regular distribution that is everywhere the span of locally Hamiltonian vector fields is automatically integrable, which we can use to define a Lie–Weinstein dual pair. The integrability of the polar distribution is a direct consequence of the regularity of the distribution that it is coming from. If the dimension of the leaves of the (generalized) integrable distribution that we start with changes from point to point —as it occurs very frequently in the Poisson symmetric case—, the associated polar distribution is in general not integrable, making it useless to define a dual pair.

In this paper we will provide a new notion of polarity and dual pair that is well defined in the absence of the regularity hypotheses needed in the classical statements. These new concepts will prove useful in recovering some of the classical results in singular situations. We will also use them to identify pseudogroups of local Poisson transformations that behave particularly well and that we will call ***von Neumann pseudogroups***. The notation has been chosen according to the resemblance of the defining properties of these pseudogroups with the von Neumann or double commutant relation for the ∗-algebras of bounded operators on a Hilbert space.

We will pay special attention to the transformation groups associated to the canonical actions of Lie groups on Poisson manifolds. More specifically, we will find various Poisson actions that are guaranteed to produce von Neumann pairs. The von Neumann character of a canonical group action has proved to be very important [O02] at the time of using the associated transformation group to implement symplectic reduction in the framework of the so called ***optimal momentum map***.

The reader will notice that in our work we will deal with various quotients that from the topological point of view may be very complex and unmanageable. Nevertheless, the use of the properties of the algebras of *smooth* functions that can be associated to these quotients allows us to extract information about this otherwise poorly behaved sets. This author is aware that the approach taken in this paper does not follow the so called *noncommutative geometry program* [Co94] that in these situations proposes the study of a $C^*$–algebra that can be associated to the equivalence relation that generates the quotient, rather than the quotient itself. This extremely suggestive approach to the problem, as well as its links to groupoid theory, will be pursued in a future work. In the same fashion in which the classical notions of polarity and Morita equivalence have been generalized to the context of $C^*$–algebras by Rieffel [Rie74, Rie74b, Land98] under the name of strong Morita equivalence, we are sure that the same can be achieved for the ideas that we present in the following pages.

The paper is structured as follows: Section 2 introduces some mathematical prerequisites that will be needed later on in the paper but, most importantly, introduces some non–standard terminology that will be strongly used in the exposition. The reader should pay special attention to Definition 2.6 and the conventions in Subsection 2.3. The main concepts of the paper are presented in Section 3. More specifically, the definitions 3.1 and 3.3 introduce the notions of polarity and dual pair, and Definition 3.9 that of von Neumann pseudogroup and von Neumann pair. Section 4 studies the correspondence between the symplectic leaves of two Poisson varieties in duality and Section 5 the relation between the duality introduced in Definition 3.3 and Howe's condition. In Section 6 we show that Hamiltonian subgroups (in the sense of Definition 2.6) are very useful at the time of constructing dual and Howe pairs. Section 7 studies the von Neumann pairs constructed using canonical Lie group actions. The main results in this section are Theorems 7.2 and 7.6 that show that proper canonical Hamiltonian Lie group actions on a Poisson manifold and compact connected Lie group symplectic actions subjected to a coisotropy condition induce von Neumann pairs. In Theorem 7.9 we obtain more von Neumann pairs out of tubewise Hamiltonian symplectic Lie group actions. Finally, the appendix in Section 8 provides a quick summary of the normal form results in [OR01] that are needed in some of the proofs in the paper.

# 2   Technical preliminaries and notation

In the following paragraphs we briefly introduce the notation and technical results that we will be using throughout the paper. The expert should be aware that some of the terminology that we use is not standard. We encourage the reader to pay special attention to Definition 2.6 and to the conventions in Subsection 2.3.

Let $M$ be a smooth manifold. A transformation group $T$ of $M$ is a subgroup of the diffeomorphisms group $\mathrm{Diff}(M)$ of $M$. The ***orbit*** $T \cdot m$ under $T$ of any element $m \in M$ is defined as the set $T \cdot m :=$ $\{F(m) \mid F \in T\}$. The relation *being in the same $T$–orbit* is an equivalence relation and induces a partition of $M$ into $T$–orbits. The space of $T$–orbits will be denoted by the quotient $M/T$.

Let $\mathrm{Diff}_L(M)$ be the monoid (set with an associative operation which contains a two–sided identity element) of ***local diffeomorphisms*** of $M$. More explicitly, the elements of $\mathrm{Diff}_L(M)$ are diffeomorphisms $F : \mathrm{Dom}(F) \subset M \to F(\mathrm{Dom}(F))$ of an open subset $\mathrm{Dom}(F) \subset M$ onto its image $F(\mathrm{Dom}(F)) \subset M$. We will denote the elements of $\mathrm{Diff}_L(M)$ as pairs $(F, \mathrm{Dom}(F))$. These local diffeomorphism can be composed using the binary operation defined as

$$(G, \mathrm{Dom}(G)) \cdot (F, \mathrm{Dom}(F)) := (G \circ F, F^{-1}(\mathrm{Dom}(G)) \cap \mathrm{Dom}(F)), \tag{2.1}$$

for all $(G, \mathrm{Dom}(G)), (F, \mathrm{Dom}(F)) \in \mathrm{Diff}_L(M)$. It is easy to see that this operation is associative and has $(\mathbb{I}, M)$, the identity map of $M$, as (unique) two sided identity element, which makes $\mathrm{Diff}_L(M)$ into a monoid. Notice that only the elements sitting in $\mathrm{Diff}(M) \subset \mathrm{Diff}_L(M)$ have an inverse since, in general, for any $(F, \mathrm{Dom}(F)) \in \mathrm{Diff}_L(M)$, we have that

$$(F^{-1}, F(\mathrm{Dom}(F))) \cdot (F, \mathrm{Dom}(F)) \quad = \quad (\mathbb{I}|_{\mathrm{Dom}(F)}, \mathrm{Dom}(F)) \tag{2.2}$$

$$(F, \mathrm{Dom}(F)) \cdot (F^{-1}, F(\mathrm{Dom}(F))) \quad = \quad (\mathbb{I}|_{F(\mathrm{Dom}(F))}, F(\mathrm{Dom}(F))). \tag{2.3}$$

Consequently, the only way to obtain the identity element $(\mathbb{I}, M)$ out of the composition of $F$ with its inverse is having $\mathrm{Dom}(F) = M$. It follows from this argument that $\mathrm{Diff}(M)$ is the biggest subgroup contained in the monoid $\mathrm{Diff}_L(M)$, with respect to the composition law (2.1). In the sequel we will frequently encounter submonoids $T_L$ of $\mathrm{Diff}_L(M)$ that contain the global identity element $(\mathbb{I}, M)$ and that satisfy the following property: for any $F : \mathrm{Dom}(F) \to F(\mathrm{Dom}(F))$ in $T_L$ there exists another element $F^{-1} : F(\mathrm{Dom}(F)) \to \mathrm{Dom}(F)$ also in $T_L$ that satisfies the identities (2.2) and (2.3). Such submonoids will be referred to as ***pseudogroups*** of $\mathrm{Diff}_L(M)$. The importance of these pseudogroups is that they have an orbit space associated. Indeed, if $T_L$ is a pseudogroup we define the orbit $T_L \cdot m$ under $T_L$ of any element $m \in M$ as the set $T_L \cdot m := \{F(m) \mid F \in T_L, m \in \mathrm{Dom}(F)\}$. $T_L$ being a pseudogroup implies that the relation *being in the same $T_L$–orbit* is an equivalence relation and induces a partition of $M$ into $T_L$–orbits. The space of $T_L$–orbits will be denoted by $M/T_L$.

**Definition 2.1** *Let $M$ be a smooth manifold and $T_L$ be one of its transformation pseudogroups. In the sequel we will use the following terminology:*

- *$T_L$ is **integrable** when its orbits are initial submanifolds of $M$, that is, if $N$ is an orbit of $T_L$ and $i : N \to M$ is the canonical injection then, for any manifold $Z$, a mapping $f : Z \to N$ is smooth iff $i \circ f : Z \to M$ is smooth.*

- *A smooth function $f \in C^\infty(M)$ is $T_L$–**invariant** when for any $(F, \mathrm{Dom}(F)) \in T_L$ we have that $f \circ F = f|_{\mathrm{Dom}(F)}$ and we denote by $C^\infty(M)^{T_L}$ the set of $T_L$–invariant functions on $M$.*

- *An open subset $U \subset M$ is said to be $T_L$–**invariant** or $T_L$–**saturated** when for any $(F, \mathrm{Dom}(F)) \in T_L$ and any $z \in \mathrm{Dom}(F) \cap U$ we have that $F(z) \in U$.*

- *The pseudogroup $T_L$ has the **extension property** when any $T_L$–invariant function $f \in C^\infty(U)^{T_L}$ defined on any $T_L$–invariant open subset $U$ has the following feature: for any $z \in U$, there is a $T_L$–invariant open neighborhood $V \subset U$ of $z$ and a $T_L$–invariant smooth function $F \in C^\infty(M)^{T_L}$ such that $f|_V = F|_V$.*

- *Finally, we say that $C^\infty(M)^{T_L}$ **separates** the $T_L$–orbits when the following condition is satisfied: if two orbits $T_L \cdot x, T_L \cdot y \in M/T_L$ are such that $f(T_L \cdot x) = f(T_L \cdot y)$ for all $f \in C^\infty(M)^{T_L}$, then $T_L \cdot x = T_L \cdot y$ necessarily.*

If $(M, \{\cdot, \cdot\})$ is a Poisson manifold, we will denote by $\mathcal{P}(M)$ the group of **Poisson automorphisms** of $M$ and by $\mathcal{P}_L(M)$ the pseudogroup of **local Poisson diffeomorphisms** of $M$. One of the main ingredients of our work in this paper will be the (finite or infinite–dimensional) subgroups and pseudosubgroups of $\mathcal{P}(M)$ and $\mathcal{P}_L(M)$, respectively, many of which will be obtained out of integrable generalized distributions on $M$. The following paragraphs review their construction.

## 2.1 Generalized distributions

We quickly review some well known facts about generalized distributions defined by families of vector fields. The standard references for this topic are [St74a, St74b], and [Su73]. We will follow the notation of [LM87].

Let $M$ be a manifold and $\mathcal{D}$ be an everywhere defined family of vector fields. By *everywhere defined* we mean that for every $m \in M$ there exists $X \in \mathcal{D}$ such that $m \in \mathrm{Dom}(X)$. The domains $\mathrm{Dom}(X) \subset M$, $X \in \mathcal{D}$, will be taken open in $M$. Let $D$ be the generalized distribution on $M$ constructed by association, to any point $z \in M$, the subspace $D(z)$ of $T_z M$ given by

$$D(z) = \mathrm{span}\{X(z) \in T_z M \,|\, X \in \mathcal{D} \text{ and } z \in \mathrm{Dom}(X)\}.$$

We will say that $D$ is the generalized distribution **spanned** by $\mathcal{D}$. Note that the dimension of $D$ may not be constant; the dimension of $D(z)$ is called the **rank** of the distribution $D$ at $z$. Any distribution defined in this way is **smooth** in the sense that for any $z \in M$ and any $v \in D(z)$ there is a smooth vector field $X$ tangent to $D$ defined in a neighborhood of $z$ and such that $X(z) = v$. An immersed connected submanifold $N$ of $M$ is said to be an **integral submanifold** of the distribution $D$ if, for every $z \in N$, $T_z i(T_z N) \subset D(z)$, where $i : N \to M$ is the canonical injection. The integral submanifold $N$ is said to be of **maximal dimension** at a point $z \in N$ if $T_z i(T_z N) = D(z)$. A **maximal integral submanifold** $N$ of $D$ is an integral manifold everywhere of maximal dimension such that any other integral submanifold of $D$, which is everywhere of maximal dimension and contains $N$, is equal to $N$. The generalized distribution $D$ is said to be **integrable** if, for every point $z \in M$, there exists a maximal integral submanifold of $D$ which contains $z$. This submanifold is usually referred to as the **leaf through** $z$ of the distribution $D$. The leaves of an integrable distribution are initial submanifolds of $M$ [Mi00].

When the distribution $D$ generated by the family of vector fields $\mathcal{D}$ is integrable, a very useful characterization of its integral manifolds can be given. In order to describe it we introduce some terminology following [LM87].

Let $X$ be a vector field defined on an open subset $\mathrm{Dom}(X)$ of $M$ and $F_t$ be its flow. For any fixed $t \in \mathbb{R}$ the domain $\mathrm{Dom}(F_t)$ of $F_t$ is an open subset of $\mathrm{Dom}(X)$ such that $F_t : \mathrm{Dom}(F_t) \to \mathrm{Dom}(F_{-t})$ is a diffeomorphism. If $Y$ is a second vector field defined on the open set $\mathrm{Dom}(Y)$ with flow $G_t$ we can consider, for two fixed values $t_1, t_2 \in \mathbb{R}$, the composition of the two diffeomorphisms $F_{t_1} \circ G_{t_2}$ as defined on the open set $\mathrm{Dom}(G_{t_2}) \cap (G_{t_2})^{-1}(\mathrm{Dom}(F_{t_1}))$ (which may be empty).

The previous prescription allows us to inductively define the composition of an arbitrary number of locally defined flows. We will obviously be interested in the flows associated to the vector fields in $\mathcal{D}$ that define the distribution $D$. The following sentences describe some important conventions that we will use all over the paper. Let $k \in \mathbb{N}$, $k > 0$, be an integer, $\mathcal{X}$ be an ordered family $\mathcal{X} = (X_1, \ldots, X_k)$ of $k$ elements of $\mathcal{D}$, and $T$ be a $k$–tuple $T = (t_1, \ldots, t_k) \in \mathbb{R}^k$ such that $F_t^i$ denotes the (locally defined) flow of $X_i$, $i \in \{1, \ldots, k\}$, $t_i$. We will denote by $\mathcal{F}_T$ the locally defined diffeomorphism $\mathcal{F}_T = F_{t_1}^1 \circ F_{t_2}^2 \circ \cdots \circ F_{t_k}^k$ constructed using the above given prescription. Any local diffeomorphism from an open subset of $M$ onto another open subset of $M$ of the form $\mathcal{F}_T$ is said to be **generated** by the family $\mathcal{D}$. It can be proven that the composition of local diffeomorphisms generated by $\mathcal{D}$ and the inverses of local diffeomorphisms generated by $\mathcal{D}$ are themselves local diffeomorphisms generated by $\mathcal{D}$ [LM87, Proposition 3.3, Appendix 3]. In other words, the local diffeomorphisms generated by $\mathcal{D}$ ***form a pseudogroup*** of $\mathrm{Diff}_L(M)$ that we will denote by $G_\mathcal{D}$. For any point $x \in M$, the symbol $G_\mathcal{D} \cdot x$ will denote the $G_\mathcal{D}$–orbit going through the point $x \in M$ and $M/G_\mathcal{D}$ the space of $G_\mathcal{D}$–orbits. In some occasions and in order to emphasize the local nature of the elements of $G_\mathcal{D} \subset \mathrm{Diff}_L(M)$ we will write them as pairs of the form $(\mathcal{F}_T, \mathrm{Dom}(\mathcal{F}_T))$.

**Theorem 2.2** *Let $D$ be a differentiable generalized distribution on the smooth manifold $M$ spanned by an everywhere defined family of vector fields $\mathcal{D}$. The following properties are equivalent:*

(i) *The distribution $D$ is invariant under the local diffeomorphisms generated by $\mathcal{D}$, that is, for each $\mathcal{F}_T \in G_\mathcal{D}$ generated by $\mathcal{D}$ and for each $z \in M$ in the domain of $\mathcal{F}_T$,*

$$T_z \mathcal{F}_T(D_z) = D_{\mathcal{F}_T(z)}.$$

(ii) *The distribution $D$ is integrable and its integral manifolds are the $G_\mathcal{D}$–orbits. Consequently, the space $M/D$ of leaves of $D$ equals $M/G_\mathcal{D}$.*

**Proof.** See [St74a, St74b, Su73]. For a compact presentation combine Theorems 3.9 and 3.10 in the Appendix 3 of [LM87]. ∎

**Notation:** In the sequel, we will use a notation consistent with the symbols just introduced: the calligraphic type $\mathcal{D}$ will denote a family of vector fields, the roman $D$ will be the associated distribution, and $G_\mathcal{D}$ will be the pseudogroup of local diffeomorphisms of $M$.

**Remark 2.3** A family $\mathcal{D}$ of locally defined vector fields on a manifold $M$ uniquely determines a pseudogroup $G_\mathcal{D}$ of local diffeomorphisms of $M$ and a generalized distribution $D$ but not the other way around, that is, a variety of families of locally defined vector fields on $M$ can be chosen in order to define the same distribution $D$. Nevertheless, if $D$ is integrable and $\mathcal{D}_1$ and $\mathcal{D}_2$ are two generating families of vector fields for $D$, the uniqueness of the maximal integral leaves of such distributions (see Theorem 2.3, page 385 of [LM87]) and the fact that by Theorem 2.2 these are given by the pseudogroup orbits, we have that for any $z \in M$, $G_{\mathcal{D}_1} \cdot z = G_{\mathcal{D}_2} \cdot z$. Consequently $M/D = M/G_{\mathcal{D}_1} = M/G_{\mathcal{D}_2}$ even though the pseudogroups $G_{\mathcal{D}_1}$ and $G_{\mathcal{D}_2}$ themselves may be different. Under some circumstances the freedom in the choice of the generating family of $D$ can be used in order to pick a family of vector fields $\mathcal{D}$ whose associated pseudogroup $G_\mathcal{D}$ is actually a subgroup of the diffeomorphism group of the manifold and hence the maximal integral manifolds of $D$ appear as group orbits. This remark motivates the introduction of the following definition. ◆

**Definition 2.4** *Let $D$ be an integrable generalized distribution on the smooth manifold $M$. We will say that $D$ is **complete** when we can choose a generating family $\mathcal{D} \in \mathfrak{X}(M)$ of $D$ made out of complete vector fields. Note that in such case the associated set of diffeomorphisms $G_\mathcal{D}$ forms a subgroup of $\mathrm{Diff}(M)$. If $\mathcal{D}$ is a generating family of $D$ that contains a subfamily that still generates $D$ and is made of complete vector fields then we say that $\mathcal{D}$ is **completable**; any such subfamily will be called a **completion** of $\mathcal{D}$.*

As we said, when we have an integrable distribution $D$ spanned by an everywhere defined family of vector fields $\mathcal{D}$ on the manifold $M$, its maximal integrable manifolds can be characterized as the orbits of the associated pseudogroup $G_\mathcal{D} \subset \mathrm{Diff}_L(M)$. This facts allows us to phrase Definition 2.1 in the context of distributions.

**Definition 2.5** *Let $D$ be an integrable generalized distribution on the smooth manifold $M$ spanned by an everywhere defined family of vector fields $\mathcal{D}$.*

- *A smooth function $f \in C^\infty(M)$ is $D$–**invariant** if it belongs to $C^\infty(M)^{G_\mathcal{D}}$. We will denote the set of $D$–invariant functions by $C^\infty(M)^D$.*

- *An open subset $U \subset M$ is called $D$–**invariant** or $D$–**saturated** if it is $G_\mathcal{D}$–invariant.*

- *The distribution $D$ has the **extension property** when the pseudogroup $G_\mathcal{D}$ has it.*

- *Finally, we say that $C^\infty(M)^D$ **separates** the integral leaves of $D$ when $C^\infty(M)^{G_\mathcal{D}}$ separates the $G_\mathcal{D}$–orbits.*

## 2.2 Subgroups and pseudosubgroups of the Poisson automorphisms group

**Definition 2.6** *Let $A$ be a subgroup of the Poisson automorphisms group $\mathcal{P}(M)$ of the Poisson manifold $(M, \{\cdot, \cdot\})$. We will denote by $C^\infty(M)^A$ the set of $A$–invariant smooth functions on $M$ and by $\left(C^\infty(M)^A\right)^c$ the centralizer of $C^\infty(M)^A$ with respect to the Poisson algebra induced by the bracket $\{\cdot, \cdot\}$ on $C^\infty(M)$.*

**(i)** *The subgroup $A$ is **strongly Hamiltonian** when every element $g \in A$ can be written as $g = F_{t_1}^1 \circ F_{t_2}^2 \circ \cdots \circ F_{t_k}^k$, with $F_{t_i}^i$ the flow of a Hamiltonian vector field $X_{h_i}$ associated to a function $h_i$ in the centralizer $\left(C^\infty(M)^A\right)^c$.*

**(ii)** *The subgroup $A$ is **weakly Hamiltonian** when for every element $g \in A$ and any $m \in M$ we can write $g \cdot m = F_{t_1}^1 \circ F_{t_2}^2 \circ \cdots \circ F_{t_k}^k(m)$, with $F_{t_i}^i$ the flow of a Hamiltonian vector field $X_{h_i}$ associated to a function $h_i \in \left(C^\infty(M)^A\right)^c$.*

**(ii)** *The subgroup $A$ is **tubewise strongly (resp. weakly) Hamiltonian** when for every element $g \in A$ and any $m \in M$ there is an $A$–invariant neighborhood $U$ of $m$ such that we can write $g = F_{t_1}^1 \circ F_{t_2}^2 \circ \cdots \circ F_{t_k}^k$ (resp. $g \cdot m = F_{t_1}^1 \circ F_{t_2}^2 \circ \cdots \circ F_{t_k}^k(m)$), with $F_{t_i}^i$ the flow of a Hamiltonian vector field $X_{h_i}$ associated to a function $h_i \in \left(C^\infty(U)^A\right)^c$.*

**Example 2.7 Connected Lie group actions with a momentum map are strongly Hamiltonian:** Let $G$ be a connected Lie group acting canonically on the Poisson manifold $(M, \{\cdot, \cdot\})$ via the map $\Phi : G \times M \to M$. The term ***canonical*** means that for any $g \in G$ and any $f, h \in C^\infty(M)$ we have that $\Phi_g^*\{f, h\} = \{\Phi_g^* f, \Phi_g^* h\}$. Suppose that the $G$-action has a momentum map $\mathbf{J} : M \to \mathfrak{g}^*$ associated. Let $A_G \subset \mathcal{P}(M)$ be the subgroup of $\mathcal{P}(M)$ defined by $A_G := \{\Phi_g : M \to M \mid g \in G\}$. Then, $A_G$ is a Hamiltonian subgroup of $\mathcal{P}(M)$. Indeed, by the connectedness of $A$, every element $g \in G$ can be written as $g = \exp \xi_1 \cdot \ldots \cdot \exp \xi_n$, with $\xi_i \in \mathfrak{g}$ in the Lie algebra $\mathfrak{g}$ of $G$. Consequently, $\Phi_g = F_1^{\xi_1} \circ F_1^{\xi_2} \circ \cdots \circ F_1^{\xi_n}$, with $F_t^{\xi_i}$ the flow of $X_{\langle \mathbf{J}, \xi_i \rangle}$. But, by Noether's Theorem $\langle \mathbf{J}, \xi_i \rangle \in \left(C^\infty(M)^G\right)^c$. ◆

**Example 2.8 A weakly and tubewise Hamiltonian group that is not Hamiltonian:** Let $M = S^1 \times S^1 = \mathbb{T}^2$ be the two torus with the symplectic form $\omega = d\theta_1 \wedge d\theta_2$ given by its area form. Let $G = S^1$ acting canonically on $M$ by $e^{i\phi} \cdot (e^{i\theta_1}, e^{i\theta_2}) := (e^{i(\phi + \theta_1)}, e^{i\theta_2})$ and $A_{S^1}$ be the associated subgroup of $\mathcal{P}(\mathbb{T}^2)$. It is easy to see that in this case, every $S^1$–invariant smooth function $f$ can be written as $f(e^{i\theta_1}, e^{i\theta_2}) = g(e^{i\theta_2})$, with $g \in C^\infty(S^1)$. Its associated Hamiltonian vector field is given by $X_f = \frac{\partial g}{\partial \theta_2} \frac{\partial}{\partial \theta_1}$. With these remarks at hand it is easy to see that $A_{S^1}$ is weakly Hamiltonian, tubewise strongly Hamiltonian, but *not* strongly Hamiltonian. ◆

**Remark 2.9** In Sections 7 and 8.3 we will study in detail some conditions under which the subgroups of the Poisson automorphism group of a manifold induced by a Lie group action are weakly and tubewise Hamiltonian. For instance, the weakly Hamiltonian character of the previous example is a corollary of one of the results that we will present (Theorem 7.6). ♦

As we already said, in many cases we will deal with pseudogroups of $\mathcal{P}_L(M)$ obtained out of integrable distributions. The following result, whose proof is a straightforward corollary of Proposition 10.3 in page 121 of [LM87], characterizes the integrable distributions whose associated pseudogroup of local diffeomorphisms lies in $\mathcal{P}_L(M)$.

**Proposition 2.10** *Let $D$ be an integrable distribution on the Poisson manifold $(M, \{\cdot, \cdot\})$ spanned by the family of vector fields $\mathcal{D}$. Then, the associated pseudogroup $G_{\mathcal{D}}$ of local diffeomorphisms of $M$ lies in $\mathcal{P}_L(M)$ iff one of the following equivalent conditions are satisfied:*

**(i)** *For any $X \in \mathcal{D}$ and any $f, g \in C^\infty(M)$, we have that*

$$X[\{f, g\}] = \{X[f], g\} + \{f, X[g]\}.$$

**(ii)** *If $B \in \Lambda^2(T^*M)$ is the Poisson tensor of $(M, \{\cdot, \cdot\})$ then, for any $X \in \mathcal{D}$ we have that $\mathcal{L}_X \Lambda = 0$. The symbol $\mathcal{L}$ denotes Lie derivation.*

The integrable distributions that fall in the category described in the previous proposition will be called ***Poisson distributions***. This denomination is sometimes used [MR86, OR98] to refer to distributions that satisfy the following condition:

If $f, g \in C^\infty(M)$ are such that $\mathbf{d}f|_D = \mathbf{d}g|_D = 0$, then $\mathbf{d}\{f, g\}|_D = 0$.

Poisson distributions always have this property but the converse is in general not true.

## 2.3 Smooth functions and Poisson structures in quotient spaces

**Definition 2.11** *A pair $(X, C^\infty(X))$, where $X$ is a topological space and $C^\infty(X) \subset C^0(X)$ is a subalgebra of the algebra of continuous functions in $X$, is called a **variety** with **smooth functions** $C^\infty(X)$. If $Y \subset X$ is a subset of $X$, the pair $(Y, C^\infty(Y))$ is said to be a **subvariety** of $(X, C^\infty(X))$, if $Y$ is a topological space endowed with the relative topology defined by that of $X$ and*

$$C^\infty(Y) = \{f \in C^0(Y) \mid f = F|_Y \text{ for some } F \in C^\infty(X)\}.$$

*Sometimes $C^\infty(Y)$ is called the set of **Whitney smooth functions** on $Y$ with respect to $X$ and is denoted by $W^\infty(Y)$. A map $\varphi : X \to Z$ between two varieties is said to be smooth when it is continuous and $\varphi^* C^\infty(Z) \subset C^\infty(X)$.*

In our discussion we will be interested in the varieties obtained as the space of orbits of the action of a pseudogroup $A$ of the local diffeomorphisms group $\text{Diff}_L(M)$ of a smooth manifold $M$; this space will be denoted by $M/A$ and we will consider it as a topological space endowed with the quotient topology. The pair $(M/A, C^\infty(M/A))$ is a variety whose algebra of smooth functions $C^\infty(M/A)$ is defined by the requirement that the canonical projection $\pi_A : M \to M/A$ is a smooth map, that is,

$$C^\infty(M/A) := \{f \in C^0(M/A) \mid f \circ \pi_A \in C^\infty(M)\}.$$

Notice that by the definition of the topology on $M/A$, the projection $\pi_A$ is continuous and, moreover, it is an open map. Indeed, if $U$ is an open set in $M$, $\pi_A(U)$ is open if and only if $\pi_A^{-1}(\pi_A(U))$ is open. Since $\pi_A^{-1}(\pi_A(U)) = A \cdot U = \cup_{\phi \in A} \phi(U \cap \mathrm{Dom}(\phi))$, $\pi_A^{-1}(\pi_A(U))$ is a union of open sets and therefore open.

In our discussion we will often work with open $A$–invariant subsets $U \subset M$ and their projections onto the orbit space $\pi_A(U) = U/A$. In principle, there are two ways to endow such an open subset of $M/A$ with a variety structure that in general do not coincide. Firstly, we can think of the variety $(U/A, C^\infty(U/A))$ with

$$C^\infty(U/A) := \{f \in C^0(U/A) \mid f \circ \pi_A|_U \in C^\infty(U)^A\}. \tag{2.4}$$

However, we can also think of $U/A$ as a subvariety of $M/A$. In that case we will denote it as $(U/A, W^\infty(U/A))$, with

$$\begin{aligned} W^\infty(U/A) &:= \{f \in C^0(U/A) \mid f = F|_{U/A} \text{ for some } F \in C^\infty(M/A)\} \\ &= \{f \in C^0(U/A) \mid f \circ \pi_A|_U = G|_U \text{ for some } G \in C^\infty(M)^A\}. \end{aligned}$$

We have the inclusion $W^\infty(U/A) \subset C^\infty(U/A)$ that in general is strict.

**Notational convention**: In all that follows and unless it is indicated otherwise we will consider the quotients of the form $U/A$ as varieties $(U/A, C^\infty(U/A))$, with $C^\infty(U/A)$ as in (2.4).

If $M$ happens to be a Poisson manifold with bracket $\{\cdot, \cdot\}$, and $A \subset \mathcal{P}_L(M)$ is a pseudogroup of $\mathcal{P}_L(M)$ then, the pair $(C^\infty(M/A), \{\cdot, \cdot\}_{M/A})$ is a well–defined Poisson algebra (also referred to as **Poisson variety**), with bracket $\{\cdot, \cdot\}_{M/A}$ given by

$$\{f, g\}_{M/A}(\pi_A(m)) = \{f \circ \pi_A, g \circ \pi_A\}(m), \tag{2.5}$$

for every $m \in M$ and any $f, g \in C^\infty(M/A)$. Analogously, if $U$ is a $A$–invariant open subset of $M$, the variety $(U/A, C^\infty(U/A))$ can be endowed with a Poisson variety structure by defining a Poisson bracket on $C^\infty(U/A)$ by restriction of that in $C^\infty(M)$, namely,

$$\{f, g\}_{U/A}(\pi_A(m)) = \{f \circ \pi_A, g \circ \pi_A\}_U(m),$$

for every $m \in U$ and any $f, g \in C^\infty(U/A)$. The symbol $\{\cdot, \cdot\}_U$ denotes the restriction of the bracket on $M$ to the open subset $U$.

The term Poisson variety is also encountered in the context of the algebraic geometric treatment of integrable systems. See for instance [Vanh96]. This concept does not in general coincide with ours.

## 3   Dual pairs

### 3.1   Polarity and dual pairs

We now introduce the notion of **polarity,** which we will use to give our definition of **dual pair.** All along this section we will be working on a smooth Poisson manifold $(M, \{\cdot, \cdot\})$.

**Definition 3.1** *Let $(M, \{\cdot, \cdot\})$ be a Poisson manifold, $A \subset \mathcal{P}_L(M)$ be a pseudosubgroup of its local Poisson diffeomorphism pseudogroup, and $(C^\infty(M/A), \{\cdot, \cdot\}_{M/A})$ be the associated quotient Poisson* *variety. Let $\mathcal{A}'$ be the set of Hamiltonian vector fields associated to all the elements of $C^\infty(U)^A$, for all the open $A$–invariant subsets $U$ of $M$, that is,*

$$\mathcal{A}' = \left\{ X_f \mid f \in C^\infty(U)^A, \text{ with } U \subset M \text{ open and } A\text{–invariant} \right\}. \tag{3.1}$$

*The distribution $A'$ associated to the family $\mathcal{A}'$ will be called the **polar distribution** defined by $A$ (or equivalently the **polar of** $A$). Any generating family of vector fields for $A'$ will be called a **polar family** of $A$. The family $\mathcal{A}'$ will be called the **standard polar family** of $A$. A pseudogroup of local Poisson diffeomorphisms associated to any polar family of $A$ will be referred to as a **polar pseudogroup** induced by $A$. The polar pseudogroup $G_{\mathcal{A}'} \subset \mathcal{P}_L(M)$ induced by the standard polar family $\mathcal{A}'$ will be called the **standard polar pseudogroup**.*

**Remark 3.2** If the pseudosubgroup $A$ has the extension property, there is a simpler polar family, we will call it $\mathcal{A}'_{ext}$, that can be used to generate $A'$, namely

$$\mathcal{A}'_{ext} = \left\{ X_f \mid f \in C^\infty(M)^A \right\}. \quad \blacklozenge$$

**Definition 3.3** *Let $(M, \{\cdot, \cdot\})$ be a Poisson manifold and $A, B \subset \mathcal{P}_L(M)$ be two pseudosubgroups of its local Poisson diffeomorphism pseudogroup. We say that the diagram*

$$
\begin{array}{ccc}
 & (M, \{\cdot, \cdot\}) & \\
{\scriptstyle \pi_A} \swarrow & & \searrow {\scriptstyle \pi_B} \\
(M/A, \{\cdot, \cdot\}_{M/A}) & & (M/B, \{\cdot, \cdot\}_{M/B})
\end{array}
$$

*is a **dual pair** on $(M, \{\cdot, \cdot\})$ when the polar distributions $A'$ and $B'$ are integrable and they satisfy that*

$$M/A' = M/B \text{ and } M/B' = M/A. \tag{3.2}$$

*The Poisson manifold $(M, \{\cdot, \cdot\})$ is called the **equivalence bimodule** of the dual pair.*

**Remark 3.4** When in (3.2) we state that $M/A' = M/B$ we mean that the partition of $M$ on $B$–orbits coincides with that of $M$ on $A'$–leaves. In general, this condition can hold without $B$ being equal to $G_{\mathcal{A}'}$ as pseudogroups; only the orbit spaces are required to be equal. Notice also that two pseudogroups $A$ and $B$ in duality are necessarily integrable. $\blacklozenge$

The following examples justify the choice of words in the previous definition.

**Example 3.5 The polar of a regular distribution and the relation with Lie's polarity.** In this example we compare the notion of polarity of Definition 3.1 with the polarity introduced by Lie [Lie90] that we described in the introduction. Let $D$ be an integrable regular distribution on the symplectic manifold $(M, \omega)$ that is the span of an everywhere defined family $\mathcal{D}$ of local Hamiltonian vector fields. As we recalled in the introduction, the space of leaves $M/D$ is a Poisson manifold and the canonical projection $\pi_D : M \to M/D$ is a Poisson surjective submersion. Lie's polar distribution $D^\omega$ is defined by

$$D^\omega(m) := \{v \in T_m M \mid \omega(m)(v, w) = 0 \text{ for all } w \in D(m)\}.$$

Since the vector fields in $\mathcal{D}$ are Hamiltonian, the associated pseudosubgroup $G_{\mathcal{D}}$ of transformations lies in $\mathcal{P}_L(M)$ and due to the integrability of $D$ we have that $M/D = M/G_{\mathcal{D}}$.

We will show that in this situation the polar $G'_{\mathcal{D}}$ of $G_{\mathcal{D}}$ coincides with $D^\omega$. We start the argument with the statement of a lemma whose proof can be found in Appendix 8.4

**Lemma 3.6** *Let $D$ be a smooth integrable regular distribution on the manifold $M$. Then, $D$ has the extension property.*

Now, since $D$ has the extension property, the polar distribution $G'_{\mathcal{D}}$ is generated by the family of globally defined vector fields (see Remark 3.2)

$$\mathcal{D}'_{ext} := \{X_f \mid f \in C^\infty(M)^D\}.$$

At the same time, since the projection $\pi_D$ is a surjective submersion we have that

$$(\ker T_m \pi_D)^\circ = \operatorname{span}\{\mathbf{d}\,(f \circ \pi_D)\,(m) \mid f \in C^\infty(M/D)\}$$

and consequently

$$D^\omega(m) = D(m)^\omega = (\ker T_m \pi_D)^\omega = \{X_{f \circ \pi_D}(m) \mid f \in C^\infty(M/D)\} = \{X_g(m) \mid g \in C^\infty(M)^D\} = G'_{\mathcal{D}}(m),$$

as required. We emphasize that, as we will see in Example 7.8, the regularity of $D$ does not imply that of its polar $D^\omega = G'_{\mathcal{D}}$. ◆

**Example 3.7 Lie–Weinstein dual pairs with connected fibers.** Let $(P_1, \{\cdot, \cdot\}_{P_1}) \overset{\pi_1}{\leftarrow} (M, \omega) \overset{\pi_2}{\rightarrow} (P_2, \{\cdot, \cdot\}_{P_2})$ be a Lie–Weinstein dual pair, that is, $\pi_1$ and $\pi_2$ are surjective Poisson submersions such that the distributions $\ker T\pi_1$ and $\ker T\pi_2$ are symplectically orthogonal. We will show that if we assume that $\pi_1$ and $\pi_2$ have connected fibers then, we can realize the diagram $(P_1, \{\cdot, \cdot\}_{P_1}) \overset{\pi_1}{\leftarrow} (M, \omega) \overset{\pi_2}{\rightarrow} (P_2, \{\cdot, \cdot\}_{P_2})$ as a dual pair in the sense of Definition 3.3. Indeed, notice that since $\pi_1$ and $\pi_2$ are surjective submersions then:

$$(\ker T_m \pi_1)^\circ \quad = \operatorname{span}\{\mathbf{d}\,(f \circ \pi_1)\,(m) \mid f \in C^\infty(P_1)\}, \tag{3.3}$$

$$(\ker T_m \pi_2)^\circ \quad = \operatorname{span}\{\mathbf{d}\,(f \circ \pi_2)\,(m) \mid f \in C^\infty(P_2)\}, \tag{3.4}$$

where the symbol $(\ker T_m \pi_1)^\circ$ denotes the annihilator of $\ker T_m \pi_1$ in $T_m^* M$. These equalities are easy to prove just by taking the local projection coordinates associated to the submersions $\pi_1$ and $\pi_2$. Now, if $B$ is the (non degenerate) Poisson tensor associated to $(M, \omega)$ and $B^\sharp : T^*M \to TM$ is the vector bundle map associated to it, we can write:

$$\ker T_m \pi_1 \quad = (\ker T_m \pi_2)^\omega = B^\sharp(m)\left((\ker T_m \pi_2)^\circ\right) = \operatorname{span}\{X_{f \circ \pi_2}(m) \mid f \in C^\infty(P_2)\} \tag{3.5}$$

$$\ker T_m \pi_2 \quad = (\ker T_m \pi_1)^\omega = B^\sharp(m)\left((\ker T_m \pi_1)^\circ\right) = \operatorname{span}\{X_{f \circ \pi_1}(m) \mid f \in C^\infty(P_1)\}. \tag{3.6}$$

Let $\mathcal{A}$ and $\mathcal{B}$ be the families of vector fields on $M$ given by

$$\mathcal{A} = \operatorname{span}\{X_{f \circ \pi_2} \mid f \in C^\infty(P_2)\} \text{ and } \mathcal{B} = \operatorname{span}\{X_{f \circ \pi_1} \mid f \in C^\infty(P_1)\},$$

and $A$ and $B$ be the associated distributions that, as a consequence of the relations (3.5) and (3.6), are guaranteed to be integrable since the level sets of $\pi_1$ and $\pi_2$ integrate them. Moreover, the connectedness hypotheses on the fibers of $\pi_1$ and $\pi_2$ allow us to make the natural identifications:

$$P_1 \simeq M/\ker T\pi_1 = M/A = M/G_{\mathcal{A}} \quad \text{and} \quad P_2 \simeq M/\ker T\pi_2 = M/B = M/G_{\mathcal{B}}.$$

Using these identifications we can rewrite the Lie–Weinstein dual pair $(P_1, \{\cdot, \cdot\}_{P_1}) \overset{\pi_1}{\leftarrow} (M, \omega) \overset{\pi_2}{\rightarrow} (P_2, \{\cdot, \cdot\}_{P_2})$ as $(M/G_{\mathcal{A}}, \{\cdot, \cdot\}_{M/G_{\mathcal{A}}}) \overset{\pi_{G_{\mathcal{A}}}}{\leftarrow} (M, \omega) \overset{\pi_{G_{\mathcal{B}}}}{\rightarrow} (M/G_{\mathcal{B}}, \{\cdot, \cdot\}_{M/G_{\mathcal{B}}})$ which, as a corollary of the previous example, is a dual pair in our sense. Indeed, since $A = \ker T\pi_1$ is a regular integrable distribution, Example 3.5 guarantees that $G'_{\mathcal{A}} = A^\omega = (\ker T\pi_1)^\omega = \ker T\pi_2 = B$, which implies that $M/G'_{\mathcal{A}} = M/B = M/G_{\mathcal{B}}$. Analogously, it can be shown that $M/G'_{\mathcal{B}} = M/G_{\mathcal{A}}$.  ◆

As we show in the next proposition the notion of polarity is particularly well behaved when it is associated to a subgroup $A$ of $\mathcal{P}(M)$ of the Poisson diffeomorphism group.

**Proposition 3.8** *Let $(M, \{\cdot, \cdot\})$ be a Poisson manifold, $A \subset \mathcal{P}(M)$ be a subgroup of its Poisson diffeomorphism group, $A'$ be the associated polar distribution, and $G_{\mathcal{A}'} \subset \mathcal{P}_L(M)$ be the standard polar pseudogroup. Then:*

**(i)** *The group $A$ commutes with the polar $G_{\mathcal{A}'}$, that is, for any $(\mathcal{F}_T, \mathrm{Dom}(\mathcal{F}_T)) \in G_{\mathcal{A}'}$ the domain $\mathrm{Dom}(\mathcal{F}_T)$ is an open $A$–invariant set and, for any $(\phi, M) \in A$ we have that $(\mathcal{F}_T \circ \phi, \mathrm{Dom}(\mathcal{F}_T)) = (\phi \circ \mathcal{F}_T, \mathrm{Dom}(\mathcal{F}_T))$*

**(ii)** *Any element $(\mathcal{F}_T, \mathrm{Dom}(\mathcal{F}_T)) \in G_{\mathcal{A}'}$ induces a local Poisson diffeomorphism $(\bar{\mathcal{F}}_T, \pi_A(\mathrm{Dom}(\mathcal{F}_T)))$ of $(M/A, \{\cdot, \cdot\}_{M/A})$, uniquely determined by the relation $\bar{\mathcal{F}}_T \circ \pi_A = \pi_A \circ \mathcal{F}_T$, that is, the standard polar pseudogroup $G_{\mathcal{A}'}$ acts canonically on $(M/A, \{\cdot, \cdot\}_{M/A})$.*

**(iii)** *The polar distribution $A'$ is Poisson and integrable. Therefore, the leaf space $M/A'$ has a natural Poisson variety $(C^\infty(M/A'), \{\cdot, \cdot\}_{M/A'})$ associated.*

**(iv)** *$A$ acts canonically on $(M/A', \{\cdot, \cdot\}_{M/A'})$. More specifically, for any $\phi \in A$, there is a Poisson diffeomorphism $\tilde{\phi}$ of $M/A'$ uniquely determined by the relation $\tilde{\phi} \circ \pi_{A'} = \pi_{A'} \circ \phi$.*

**Proof.**   **(i)** Let $(\phi, M) \in A$, and $(\mathcal{F}_T, \mathrm{Dom}(\mathcal{F}_T)) \in G_{\mathcal{A}'}$. For the sake of simplicity in the presentation we will take $(\mathcal{F}_T, \mathrm{Dom}(\mathcal{F}_T))$ to be the flow $(F_t, \mathrm{Dom}(F_t))$ of $X_h$, with $h \in C^\infty(U)^A$ and $U$ and open $A$–invariant subset of $M$. Using the $A$–invariance of $h$ and the Poisson character of $\phi$ it is easy to see that

$$X_h = T\phi \circ X_h \circ \phi^{-1}. \tag{3.7}$$

Let now $G_t : \phi(\mathrm{Dom}(F_t)) \to \phi(F_t(\mathrm{Dom}(F_t)))$ be the local diffeomorphism defined by $G_t := \phi \circ F_t \circ \phi^{-1}$. The chain rule and expression (3.7) show that for any $z \in \mathrm{Dom}(F_t)$

$$\frac{d}{dt} G_t(\phi(z)) = (T\phi \circ X_h)(F_t(\phi^{-1}(\phi(z)))) = (T\phi \circ X_h \circ \phi^{-1})(G_t(\phi(z)) = X_h(G_t(\phi(z)). \tag{3.8}$$

The uniqueness of the flow of a vector field implies that $\phi(\mathrm{Dom}(F_t)) \subset \mathrm{Dom}(F_t)$. Since $\phi \in A$ is arbitrary, we also have that $\phi^{-1}(\mathrm{Dom}(F_t)) \subset \mathrm{Dom}(F_t)$ and, consequently $\phi(\mathrm{Dom}(F_t)) = \mathrm{Dom}(F_t)$. Expression (3.8) also implies that $F_t = G_t = \phi \circ F_t \circ \phi^{-1}$ which guarantees the commutation relation in the statement.

**(ii)** Given $(\mathcal{F}_T, \mathrm{Dom}(\mathcal{F}_T)) \in G_{\mathcal{A}'}$, the existence of the well defined map $(\bar{\mathcal{F}}_T, \pi_A(\mathrm{Dom}(\mathcal{F}_T))) = (\bar{\mathcal{F}}_T, \mathrm{Dom}(\mathcal{F}_T)/A)$ that satisfies $\bar{\mathcal{F}}_T \circ \pi_A = \pi_A \circ \mathcal{F}_T$ is guaranteed by **(i)**. Since $\mathcal{F}_T$ is a local diffeomorphism of $M$ and the projection $\pi_A$ is open and continuous, $\bar{\mathcal{F}}_T$ is necessarily continuous. We also have that $\bar{\mathcal{F}}_T^* C^\infty(\mathcal{F}_T(\mathrm{Dom}(\mathcal{F}_T))/A) \subset C^\infty(\mathrm{Dom}(\mathcal{F}_T)/A)$ since for any $f \in C^\infty(\mathcal{F}_T(\mathrm{Dom}(\mathcal{F}_T))/A)$ the map $f \circ \bar{\mathcal{F}}_T \circ \pi_A|_{\mathrm{Dom}(\mathcal{F}_T)} = f \circ \pi_A \circ \mathcal{F}_T|_{\mathrm{Dom}(\mathcal{F}_T)} \in C^\infty(\mathrm{Dom}(\mathcal{F}_T))^A$ and hence $f \circ \bar{\mathcal{F}}_T \in C^\infty(\mathrm{Dom}(\mathcal{F}_T)/A)$. Since we could do the same with $\bar{\mathcal{F}}_T^{-1}$, we conclude that the map $\bar{\mathcal{F}}_T$ is a local diffeomorphism. A straightforward verification shows that $\bar{\mathcal{F}}_T$ is also a Poisson map between $(\mathrm{Dom}(\mathcal{F}_T)/A, \{\cdot, \cdot\}_{\mathrm{Dom}(\mathcal{F}_T)/A})$ and $(\mathcal{F}_T(\mathrm{Dom}(\mathcal{F}_T))/A, \{\cdot, \cdot\}_{\mathcal{F}_T(\mathrm{Dom}(\mathcal{F}_T))/A})$.

**(iii)** First of all notice that the elements of $G_{\mathcal{A}'}$ are finite compositions of Hamiltonian flows and therefore are local Poisson diffeomorphisms. This makes of $A'$ a Poisson distribution. As to its integrability, according to Theorem 2.2 we have to show that for any $(\mathcal{F}_T, \mathrm{Dom}(\mathcal{F}_T)) \in G_{\mathcal{A}'}$ and any $m \in \mathrm{Dom}(\mathcal{F}_T)$ we have that $T_m \mathcal{F}_T(A'(m)) = A'(\mathcal{F}_T(m))$. In order to establish this equality we take $h \in C^\infty(U)^A$ with $U$ an open $A$–invariant subset of $M$. Let $V := U \cap \mathrm{Dom}(\mathcal{F}_T)$ be such that $m \in V$ and define $\mathcal{F}_T^V := \mathcal{F}_T|_V : V \to \mathcal{F}_T(V)$ and $h^V := h|_V$. Given that $V$ is an $A$–invariant open subset of $M$ and that $\mathcal{F}_T^V$ is a Poisson map we can write

$$T_m \mathcal{F}_T \left( X_h(m) \right) = T_m \mathcal{F}_T^V \left( X_{h^V}(m) \right) = T_m \mathcal{F}_T^V \left( X_{h^V \circ (\mathcal{F}_T^V)^{-1} \circ \mathcal{F}_T^V}(m) \right) = X_{h^V \circ (\mathcal{F}_T^V)^{-1}}(\mathcal{F}_T^V(m))$$

which belongs to $A'(\mathcal{F}_T(m))$ since by point **(i)**, $h^V \circ (\mathcal{F}_T^V)^{-1}$ belongs to $C^\infty(\mathcal{F}_T(V))^A$. Consequently, $T_m \mathcal{F}_T(A'(m)) \subset A'(\mathcal{F}_T(m))$.

Conversely, let $f \in C^\infty(W)^A$ be such that $\mathcal{F}_T(m) \in W$. Define $S := \mathcal{F}_T(\mathrm{Dom}(\mathcal{F}_T)) \cap W$, $f^S := f|_S$, $\mathcal{F}_T^S := \mathcal{F}_T|_{\mathcal{F}^{-1}(S)}$, then

$$X_f(\mathcal{F}_T(m)) = X_{f^S}(\mathcal{F}_T^S(m)) = T_m \mathcal{F}_T^S(X_{f^S \circ \mathcal{F}_T^S}(m)),$$

which belongs to $T_m \mathcal{F}_T(A'(m))$, as required.

**(iv)** It is a straightforward consequence of the fact, proved in **(i)**, that $A$ and $G_{\mathcal{A}'}$ commute.  ∎

## 3.2 Von Neumann dual pairs

**Definition 3.9** *Let $(M, \{\cdot, \cdot\})$ be a smooth Poisson manifold and $A \subset \mathcal{P}_L(M)$ be a pseudosubgroup of its local Poisson diffeomorphism pseudogroup $\mathcal{P}_L(M)$. We say that $A$ is **von Neumann** when the diagram $(M/A, \{\cdot, \cdot\}_{M/A}) \xleftarrow{\pi_A} (M, \{\cdot, \cdot\}) \xrightarrow{\pi_{A'}} (M/A', \{\cdot, \cdot\}_{M/A'})$ is a dual pair or, equivalently, when the distributions $A'$ and $(A')'$ are integrable and*

$$M/(A')' = M/A. \tag{3.9}$$

*In the presence of this condition we also talk of the **von Neumann pair** associated to $A \subset \mathcal{P}_L(M)$.*

**Remark 3.10** When in the previous definition $A$ is actually a subgroup of $\mathcal{P}(M)$, Proposition 3.8 automatically guarantees the integrability of $A'$.  ♦

**Remark 3.11** Von Neumann groups have connected and path connected orbits since the relation (3.9) implies that for any point $m \in M$, the orbit $A \cdot m$ coincides with $G_{\mathcal{A}''} \cdot m$ which is a connected and path connected set.  ♦

**Remark 3.12** The terminology in the previous definition has been chosen according to the similarity of condition (3.9) with the von Neumann or double commutant relation for $*$–algebras of bounded operators on a Hilbert space.  ♦

**Example 3.13 Lie group canonical actions and the optimal momentum map.** Let $G$ be a connected Lie group acting on the symplectic manifold $(M, \omega)$ in a free, proper, and canonical fashion via the map $\Phi : G \times M \to M$. The term ***canonical*** means that for any $g \in G$ we have that $\Phi_g^* \omega = \omega$. We will for the time being also assume that it has an associated equivariant momentum map that we will denote by $\mathbf{J} : M \to \mathfrak{g}^*$ whose level sets are connected. The symbol $A_G \subset \mathcal{P}(M)$ will denote the subgroup of $\mathcal{P}(M)$ defined by $A_G := \{\Phi_g : M \to M \mid g \in G\}$ and $\pi : M \to M/G$ the surjective submersion obtained by projecting $M$ onto the orbit space $M/G = M/A_G$. A calculation that is left to the reader as an exercise (see [OR02]) shows that in the presence of our hypotheses (free and proper canonical action with a momentum map) the polar distribution $A_G'$ of $A_G$ is given by

$$A_G'(m) = \ker T_m \mathbf{J} \quad \text{for all} \quad m \in M.$$

Consequently, in this particular example, the diagram $(M/A_G, \{\cdot, \cdot\}_{M/A_G}) \xleftarrow{\pi_{A_G}} (M, \omega) \xrightarrow{\pi_{A_G'}} (M/A_G', \{\cdot, \cdot\}_{M/A_G'})$ can be identified with $(M/G, \{\cdot, \cdot\}_{M/G}) \xleftarrow{\pi} (M, \omega) \xrightarrow{\mathbf{J}} (\mathfrak{g}_{\mathbf{J}}^*, \{\cdot, \cdot\}_{\mathfrak{g}_{\mathbf{J}}^*})$, where $\mathfrak{g}_{\mathbf{J}}^* =: \mathbf{J}(M)$ and $\{\cdot, \cdot\}_{\mathfrak{g}_{\mathbf{J}}^*}$ is the restriction to $\mathfrak{g}_{\mathbf{J}}^*$ of the Lie–Poisson structure of $\mathfrak{g}^*$. Given that for any $m \in M$ we have that $(\ker T_m \mathbf{J})^\omega = \mathfrak{g} \cdot m$, this diagram is a Lie–Weinstein pair with connected fibers and, by Example 3.7, a dual pair in our sense. We have thereby shown that the subgroup $A_G$ associated to a free canonical connected Lie group action that has a momentum map associated with connected fibers is von Neumann.

One of the main goals of Section 7 will be the study of $A_G$ in more realistic situations, namely when the $G$–action is not free anymore, as well as the search for situations in which the diagram $(M/A_G, \{\cdot, \cdot\}_{M/A_G}) \xleftarrow{\pi_{A_G}} (M, \{\cdot, \cdot\}) \xrightarrow{\pi_{A_G'}} (M/A_G', \{\cdot, \cdot\}_{M/A_G'})$ is still a dual pair despite the singularities in the problem. Recall that no matter how singular the $G$–action is, the right hand side leg of the previous diagram is always well defined since by Proposition 3.8 the distribution $A_G'$ is always integrable. The projection $\pi_{A_G'} : M \to M/A_G'$ will be referred to as the ***optimal momentum map*** associated to the canonical $G$–action on $M$ and will be denoted by $\mathcal{J} : M \to M/A_G'$. The reader is encouraged to check with [OR02] for a detailed study of this object. ◆

**Remark 3.14** The previous example describes a situation where it is very easy to compare the notion of polarity used by Lie and Weinstein [Lie90, W83] with ours. Suppose that the manifold $M$ is symplectic with form $\omega$ and, for the sake of simplicity, the Lie group $G$ is connected, acts properly on $M$, and has Lie algebra $\mathfrak{g}$. Let $\mathcal{D}_G$ be the family of vector fields defined by the infinitesimal generators of the $G$–action on $M$, $D_G$ be the associated distribution, and $A_G \subset \mathcal{P}(M)$ the corresponding Poisson diffeomorphism group. By definition, for any $m \in M$ we have that

$$D_G(m) = \operatorname{span}\{\xi_M(m) \mid \xi \in \mathfrak{g}\} = T_m(G \cdot m) =: \mathfrak{g} \cdot m,$$

and, by the connectedness of $G$ one has that $M/G = M/D_G = M/A_G$. If we use the classical definition [Lie90, W83], the polar of $A_G$, we will call it $A_G^\perp$, is the distribution

$$A_G^\perp(m) := (\mathfrak{g} \cdot m)^\omega = B^\sharp(m)((\mathfrak{g} \cdot m)^\circ),$$

which in general is not integrable.

We now compute the polar $A_G'$ of $A_G$ according to our definition, using the fact, whose proof can be found in [O98, OR02], that for any point $m \in M$ with isotropy subgroup $H := G_m$ we have that:

$$((T_m(G \cdot m))^\circ)^H = \operatorname{span}\{\mathbf{d}f(m) \mid f \in C^\infty(M)^G\}.$$

The symbol $((T_m(G \cdot m))^\circ)^H$ denotes the set of fixed points by the action of $H$ in the vector space $(T_m(G \cdot m))^\circ$ or, more explicitly:

$$((T_m(G \cdot m))^\circ)^H = \{v \in (T_m(G \cdot m))^\circ \mid h \cdot v = v \quad \text{for all} \quad h \in H\}.$$

By definition $A'_G$ is the distribution associated to the family of vector fields

$$\mathcal{A}'_G(m) := \{X_f(m) \mid f \in (C^\infty(U))^G, \text{ with } U \subset M \text{ open, } G\text{–invariant, } m \in U\}.$$

As the $G$–action is proper, a standard result (see [ACG91]) guarantees that any $G$–invariant function defined on a $G$–invariant open subset of $M$ admits an extension to a $G$–invariant function on $M$. This circumstance allows to simplify the definition of $A'_G$ as follows

$$
\begin{aligned}
A'_G(m): \quad &= \quad \text{span}\{X_f(m) \mid f \in (C^\infty(M))^G\} \\
&= \quad B^\sharp(m)(\{\mathbf{d}f(m) \mid f \in C^\infty(M)^G\}) = B^\sharp(m)(((\mathfrak{g} \cdot m)^\circ)^H).
\end{aligned}
$$

This distribution is always integrable. Notice that in the presence of symmetric points (that is $H \neq \{e\}$) the distributions $A_G^{\perp}$ and $A'_G$ are in general different, making the two notions of polarity not to be the same. We emphasize that even though when the situation is regular (meaning that the $G$–action is free) both notions coincide, in the singular case, the notion of polarity given in Definition 3.1 is preferable since it produces integrable distributions that can be used to define dual pairs. $\blacklozenge$

# 4 Dual pairs and symplectic leaf correspondence

It is a well known fact that the symplectic leaves of the two Poisson manifolds in the legs of a Lie–Weinstein dual pair with connected fibers are in bijection. This result was introduced in [W83]. See also the Appendix E of [Bl01] for a fully detailed proof.

In this section we will see that the situation is analogous for the two Poisson varieties in the legs of the dual pairs introduced in Definition 3.3. Nevertheless, since in this context there is no Symplectic Stratification Theorem we need to start by defining what we mean by the symplectic leaves of a quotient Poisson variety.

Let $(M, \{\cdot, \cdot\})$ be a smooth Poisson manifold, $A \subset \mathcal{P}(M)$ be a subgroup of its Poisson diffeomorphism group, and $(M/A, \{\cdot, \cdot\}_{M/A})$ be the associated quotient Poisson variety. Let $V \subset M/A$ be an open subset of $M/A$ and $h \in C^\infty(V)$ be a smooth function defined on it. If we call $U = \pi_A^{-1}(V)$ then, the vector field $X_{h \circ \pi_A|_U}$ belongs to $\mathcal{A}'$ and, by part **(ii)** of Proposition 3.8, its flow $(F_t, \text{Dom}(F_t))$ uniquely determines a local Poisson diffeomorphism $(\bar{F}_t, \pi_A(\text{Dom}(F_t)))$ of $M/A$. We will say that $(\bar{F}_t, \pi_A(\text{Dom}(F_t)))$ is the **Hamiltonian flow** associated to $h$. The symplectic leaves of $M/A$ will be defined as the reachable sets in this quotient by finite compositions of Hamiltonian flows. Since it is not clear how to defined these flows by projection of $A$–equivariant flows when $A$ is a pseudogroup of local transformations in $\mathcal{P}_L(M)$, we will restrict in this section to the case $A \subset \mathcal{P}(M)$.

**Definition 4.1** *Let* $(M, \{\cdot, \cdot\})$ *be a smooth Poisson manifold,* $A \subset \mathcal{P}(M)$ *be a subgroup of its Poisson diffeomorphism group, and* $(M/A, \{\cdot, \cdot\}_{M/A})$ *be the associated quotient Poisson variety. Given a point* $[m]_A \in M/A$, *the* **symplectic leaf** $\mathcal{L}_{[m]_A}$ *going through it is defined as the (path connected) set formed by all the points that can be reached from* $[m]_A$ *by applying to it a finite number of Hamiltonian flows associated to functions in* $C^\infty(V)$, *with* $V \subset M/A$ *any open subset of* $M/A$, *that is,*

$$\mathcal{L}_{[m]_A} := \{F_{t_1}^1 \circ F_{t_2}^2 \circ \cdots \circ F_{t_k}^k([m]_A) \mid k \in \mathbb{N}, \, F_{t_i} \text{ flow of some } X_{h_i}, \, h_i \in C^\infty(V), \, V \subset M/A \text{ open}\}.$$

*The relation* being in the same symplectic leaf *determines an equivalence relation in* $M/A$ *whose corresponding space of equivalence classes will be denoted by* $(M/A)/\{\cdot, \cdot\}_{M/A}$.

**Remark 4.2** In the paragraph preceding the previous definition the choice of the term *Hamiltonian flow* for $(\bar{F}_t, \pi_A(\mathrm{Dom}(F_t)))$ is justified by the fact that for any other function $f \in C^\infty(V)$ and any $[m]_A \in \pi_A(\mathrm{Dom}(F_t)))$ we have that

$$\frac{d}{dt} f(\bar{F}_t([m]_A)) = \{f, h\}_V(\bar{F}_t([m]_A)),\tag{4.1}$$

where $\{\cdot, \cdot\}_V$ denotes the restriction to $V$ of the Poisson bracket $\{\cdot, \cdot\}_{M/A}$. Nevertheless, expression (4.1) does not fully characterize, in general, the flow $\bar{F}_t$ since there could be other mappings for which such equality holds. This could be rephrased by saying that in the category in which we are working any function has a Hamiltonian flow associated but, unlike the smooth Poisson category, its uniqueness is not guaranteed. One result in this direction that can be easily proven by mimicking the results in page 389 of [SL91] says that if the functions in $C^\infty(M/A)$ separate the points of $M/A$ (that is, if $f(x) = f(y)$ for all $C^\infty(M/A)$, then $x = y$) then any Hamiltonian function has a unique flow satisfying the relation (4.1). ♦

Even though in Definition 4.1 we called $\mathcal{L}_{[m]_A}$ a symplectic leaf, there is in general no natural way to define on this set a smooth structure and a symplectic form that would make it a symplectic manifold. Nevertheless, there is still something we can do to justify our notation. Indeed, if we consider the set $\mathcal{L}_{[m]_A}$ as a subvariety of $M/A$ in the sense of Definition 2.11, the corresponding ring of smooth functions $C^\infty(\mathcal{L}_{[m]_A})$ given by

$$C^\infty(\mathcal{L}_{[m]_A}) = \left\{ f \in C^0(M/A) \,\middle|\, f = F|_{\mathcal{L}_{[m]_A}}, F \in C^\infty(M/A) \right\}.$$

can be endowed with a natural Poisson algebra structure via a bracket $\{\cdot, \cdot\}_{\mathcal{L}_{[m]_A}}$ that we will describe later on. It turns out that if $A$ has the extension property, the Poisson algebra $(C^\infty(\mathcal{L}_{[m]_A}), \{\cdot, \cdot\}_{\mathcal{L}_{[m]_A}})$ is non degenerate, which is the closest that we can get to being symplectic in this category; on other words, if we look at a smooth symplectic manifold $(M, \omega)$ from the point of view of its smooth functions and the Poisson bracket $\{\cdot, \cdot\}$ defined on them via the symplectic form, the symplecticity is reflected in the non degeneracy of the algebra $(C^\infty(M), \{\cdot, \cdot\})$, which is exactly the property that we will prove for $(C^\infty(\mathcal{L}_{[m]_A}), \{\cdot, \cdot\}_{\mathcal{L}_{[m]_A}})$. We make these claims more explicit in the statement of the following proposition.

**Proposition 4.3** *Let $(M, \{\cdot, \cdot\})$ be a smooth Poisson manifold, $A \subset \mathcal{P}(M)$ be a subgroup of its Poisson diffeomorphism group, and $(M/A, \{\cdot, \cdot\}_{M/A})$ be the associated quotient Poisson variety. Let $[m]_A \in M/A$ and $\mathcal{L}_{[m]_A}$ be the symplectic leaf through it. Then, the ring $C^\infty(\mathcal{L}_{[m]_A})$ can be endowed with a natural Poisson algebra structure $(C^\infty(\mathcal{L}_{[m]_A}), \{\cdot, \cdot\}_{\mathcal{L}_{[m]_A}})$ with $\{\cdot, \cdot\}_{\mathcal{L}_{[m]_A}}$ the bracket defined by*

$$\{f, g\}_{\mathcal{L}_{[m]_A}}([z]_A) := \{F, G\}_{M/A}([z]_A) = \{F \circ \pi_A, G \circ \pi_A\}(z),\tag{4.2}$$

*for any $[z]_A \in \mathcal{L}_{[m]_A}$, $f, g \in C^\infty(\mathcal{L}_{[m]_A})$, and any $F, G \in C^\infty(M/A)$ such that $F|_{\mathcal{L}_{[m]_A}} = f$ and $G|_{\mathcal{L}_{[m]_A}} = g$.*

*Moreover, if $A$ has the extension property, then the Poisson algebra $(C^\infty(\mathcal{L}_{[m]_A}), \{\cdot, \cdot\}_{\mathcal{L}_{[m]_A}})$ is non degenerate, that is, if $f \in C^\infty(\mathcal{L}_{[m]_A})$ is such that $\{f, g\}_{\mathcal{L}_{[m]_A}} = 0$ for all $g \in C^\infty(\mathcal{L}_{[m]_A})$, then $f$ is a constant function.*

**Proof.** In order to establish the first part of the Proposition it suffices to show that the bracket (4.2) is well defined or, more explicitly that its value does not depend on the extensions $F, G \in C^\infty(M/A)$ that come into its definition. Let $G' \in C^\infty(M/A)$ be another extension of $g \in C^\infty(\mathcal{L}_{[m]_A})$, $(F_t, \mathrm{Dom}(F_t))$ be the flow of the Hamiltonian vector field $X_{F \circ \pi_A}$ and $(\bar{F}_t, \pi_A(\mathrm{Dom}(F_t)))$ be the local Poisson diffeomorphism of $(M/A, \{\cdot, \cdot\}_{M/A})$, uniquely determined by the relation $\bar{F}_t \circ \pi_A = \pi_A \circ F_t$. Then,

$$
\begin{aligned}
\{F, G'\}_{M/A}([z]_A) &= \{F \circ \pi_A, G' \circ \pi_A\}(z) = -\mathbf{d}(G' \circ \pi_A)(z) \cdot X_{F \circ \pi_A}(z) \\
&= -\left.\frac{d}{dt}\right|_{t=0} G' \circ \pi_A(F_t(z)) = \left.\frac{d}{dt}\right|_{t=0} G' \circ \bar{F}_t([z]_A) \\
&= -\left.\frac{d}{dt}\right|_{t=0} G \circ \bar{F}_t([z]_A) = \{F, G\}_{M/A}([z]_A),
\end{aligned}
$$

where $G' \circ \bar{F}_t([z]_A) = G \circ \bar{F}_t([z]_A)$ because $\bar{F}_t([z]_A) \in \mathcal{L}_{[m]_A}$. Analogously, if we take another extension $F'$ of $f$ we have that for any $[z]_A \in \mathcal{L}_{[m]_A}$, $\{F', G'\}_{M/A}([z]_A) = \{F, G\}_{M/A}([z]_A)$, which proves that the bracket $\{\cdot, \cdot\}_{\mathcal{L}_{[m]_A}}$ is well defined. The rest of the defining properties of a Poisson bracket are a straightforward verification.

We now assume that $A$ has the extension property and show that the bracket $\{\cdot, \cdot\}_{\mathcal{L}_{[m]_A}}$ is non degenerate. Let $f \in C^\infty(\mathcal{L}_{[m]_A})$ be such that $\{f, g\}_{\mathcal{L}_{[m]_A}} = 0$ for all $g \in C^\infty(\mathcal{L}_{[m]_A})$. Before we proceed, the reader should notice that the definition of symplectic leaf and part **(ii)** of Proposition 3.8 imply that $\mathcal{L}_{[m]_A} = G_{\mathcal{A}'} \cdot [m]_A$. Now, since by hypothesis $A$ has the extension property, any element $\mathcal{F}_T \in G_{\mathcal{A}'}$ can be written as a finite composition of Hamiltonian flows associated to functions in $C^\infty(M)^A$ (see Remark 3.2); for the sake of simplicity we take $\mathcal{F}_T = F_t$, with $(F_t, \mathrm{Dom}(F_t))$ the flow of $X_g$, $g \in C^\infty(M)^A$. Let $G \in C^\infty(M/A)$ be the function uniquely determined by the equality $g = G \circ \pi_A$, $(\bar{F}_t, \pi_A(\mathrm{Dom}(F_t)))$ be the unique local Poisson diffeomorphism of $M/A$ defined by the relation $\bar{F}_t \circ \pi_A|_{\mathrm{Dom}(F_t)} = \pi_A \circ F_t$, and $g' = G|_{\mathcal{L}_{[m]_A}} \in C^\infty(\mathcal{L}_{[m]_A})$. We now take $F \in C^\infty(M/A)$ such that $F|_{\mathcal{L}_{[m]_A}} = f$. With all these ingredients we have that

$$
\begin{aligned}
\frac{d}{dt} f \circ \bar{F}_t([m]_A) &= \frac{d}{dt} F(\bar{F}_t([m]_A)) = \frac{d}{dt} F \circ \pi_A \circ F_t(m) = \mathbf{d}(F \circ \pi_A)(F_t(m)) \cdot X_g(F_t(m)) \\
&= \{F \circ \pi_A, g\}(F_t(m)) = \{F \circ \pi_A, G \circ \pi_A\}(F_t(m)) = \{F, G\}_{M/A}(\bar{F}_t([m]_A)) \\
&= \{f, g'\}_{\mathcal{L}_{[m]_A}}(\bar{F}_t([m]_A)) = 0,
\end{aligned}
$$

which implies that $f \circ \bar{F}_t([m]_A) = f([m]_A)$. The arbitrary character of $\bar{F}_t$ allows us to write

$$
f\left(G_{\mathcal{A}'} \cdot [m]_A\right) = f(\mathcal{L}_{[m]_A}) = f([m]_A),
$$

and thereby, the function $f \in C^\infty\left(\mathcal{L}_{[m]_A}\right)$ is constant, as required. ∎

**Remark 4.4** If in the previous Proposition we drop the hypothesis on the extension property of $A$ then the Poisson algebra $\left(C^\infty(\mathcal{L}_{[m]_A}), \{\cdot, \cdot\}_{\mathcal{L}_{[m]_A}}\right)$ is still non degenerate in the following generalized sense: if $f \in C^\infty(\mathcal{L}_{[m]_A})$ is such that for any open subset $U \subset \mathcal{L}_{[m]_A}$ and any $g \in C^\infty(U)$ we have that $\{f|_U, g\}_U = 0$, with $\{\cdot, \cdot\}_U$ the restriction to $U$ of the bracket $\{\cdot, \cdot\}_{\mathcal{L}_{[m]_A}}$, then $f$ is a constant function. ♦

**Theorem 4.5 (Symplectic leaves correspondence)** *Let $(M, \{\cdot, \cdot\})$ be a smooth Poisson manifold, $A, B \subset \mathcal{P}(M)$ be two subgroups of its Poisson diffeomorphism group, and $G_{\mathcal{A}'}, G_{\mathcal{B}'} \subset \mathcal{P}_L(M)$ be the standard polar pseudogroups. If we denote by $(M/A)/\{\cdot, \cdot\}_{M/A}$ and $(M/B)/\{\cdot, \cdot\}_{M/B}$ the space of symplectic leaves of the Poisson varieties $(M/A, \{\cdot, \cdot\}_{M/A})$ and $(M/B, \{\cdot, \cdot\}_{M/B})$, respectively, we have that:*

**(i)** *The symplectic leaves of $M/A$ and $M/B$ are given by the orbits of the $G_{\mathcal{A}'}$ and $G_{\mathcal{B}'}$ actions on $M/A$ and $M/B$, respectively, as defined in Proposition 3.8. As a consequence of this statement, we can write that*

$$(M/A)/\{\cdot, \cdot\}_{M/A} = (M/A)/G_{\mathcal{A}'} \ \text{and} \ (M/B)/\{\cdot, \cdot\}_{M/B} = (M/B)/G_{\mathcal{B}'}. \tag{4.3}$$

**(ii)** *If the diagram $(M/A, \{\cdot, \cdot\}_{M/A}) \xleftarrow{\pi_A} (M, \{\cdot, \cdot\}) \xrightarrow{\pi_B} (M/B, \{\cdot, \cdot\}_{M/B})$ is a dual pair then the map*

$$\begin{array}{rcl}
(M/A)/\{\cdot, \cdot\}_{M/A} & \longrightarrow & (M/B)/\{\cdot, \cdot\}_{M/B} \\
\mathcal{L}_{[m]_A} & \longmapsto & \mathcal{L}_{[m]_B}
\end{array} \tag{4.4}$$

*is a bijection. The symbols $\mathcal{L}_{[m]_A}$ and $\mathcal{L}_{[m]_B}$ denote the symplectic leaves in $M/A$ and $M/B$, respectively, going through the point $[m]_A$ and $[m]_B$.*

**Proof.** **(i)** It is a straightforward consequence of the definition of symplectic leaf and of the actions, spelled out in Proposition 3.8, of the standard polar pseudogroups on the quotients.

**(ii)** Given that by Proposition 3.8 $A$ and $G_{\mathcal{A}'}$ (resp. $B$ and $G_{\mathcal{B}'}$) commute, and using the duality hypothesis, we can write

$$(M/A)/\{\cdot, \cdot\}_{M/A} = (M/A)/G_{\mathcal{A}'} \simeq (M/G_{\mathcal{A}'})/A = (M/B)/A, \tag{4.5}$$

and the same relation for the subgroup $B$, that is,

$$(M/B)/\{\cdot, \cdot\}_{M/B} = (M/B)/G_{\mathcal{B}'} \simeq (M/G_{\mathcal{B}'})/B = (M/A)/B. \tag{4.6}$$

In the previous expressions $(M/B)/A$ and $(M/A)/B$ should be understood as the orbit spaces of the $A$ and $B$ actions on $M/B$ and $M/A$, respectively, inherited from considering these quotients as $M/G_{\mathcal{A}'}$ and $M/G_{\mathcal{B}'}$. More explicitly, for any $a \in A$ and any $[m]_B \in M/B$ we define $a \cdot [m]_B := a \cdot [m]_{G_{\mathcal{A}'}} = [a \cdot m]_{G_{\mathcal{A}'}} = [a \cdot m]_B$. Analogously, for any $b \in B$ and any $[m]_A \in M/A$, we define $b \cdot [m]_A := [b \cdot m]_A$. With these conventions and in view of (4.5) and (4.6) the bijective character of the map in the statement will be proved if we show that the map

$$\begin{array}{rcl}
F: \ (M/B)/A & \longrightarrow & (M/A)/B \\
[[m]_B]_A & \longmapsto & [[m]_A]_B
\end{array}$$

is a well defined bijection. It is indeed so since if $[[m]_B]_A = [[m']_B]_A$, there exists elements $a \in A$ and $b \in B$ such that $a \cdot m = b \cdot m'$ and hence $F([[m]_B]_A) = [[m]_A]_B = [[a \cdot m]_A]_B = [[b \cdot m']_A]_B = [b \cdot [m']_A]_B = [[m']_A]_B = F([[m']_B]_A)$, which shows that the map $F$ is well defined. Analogously one shows that $F$ is one to one and onto, as required. $\blacksquare$

**Remark 4.6** As a consequence of the previous theorem we can conclude that the symplectic leaves of two Poisson manifolds in the legs of a Lie–Weinstein dual pair $(P_1, \{\cdot, \cdot\}_{P_1}) \xleftarrow{\pi_1} (M, \omega) \xrightarrow{\pi_2} (P_2, \{\cdot, \cdot\}_{P_2})$ in which the projections $\pi_1$ and $\pi_2$ are complete and have connected fibers are in bijection. We say that a smooth map $T: P \to Q$ between two Poisson manifolds $P$ and $Q$ is ***complete*** if for any function $f \in C^\infty(Q)$ whose Hamiltonian vector field $X_f$ associated is complete, the vector field $X_{f \circ T} \in \mathfrak{X}(P)$ is also complete. In our context this condition shows up when we put the Lie–Weinstein dual pair in our language by making the identifications $P_1 \simeq M/G_{\mathcal{A}_c}$ and $P_2 \simeq M/G_{\mathcal{B}_c}$, with

$$\mathcal{A}_c = \text{span}\left\{X_{f \circ \pi_2} \mid f \in C_c^\infty(P_2)\right\} \ \text{and} \ \mathcal{B}_c = \text{span}\left\{X_{f \circ \pi_1} \mid f \in C_c^\infty(P_1)\right\}.$$

The subscript $c$ in $C_c^\infty(P_1)$ and $C_c^\infty(P_2)$ denotes compactly supported functions. The completeness of the projections $\pi_1$ and $\pi_2$ ensures that $G_{\mathcal{A}_c}$ and $G_{\mathcal{B}_c}$ are subgroups of $\mathcal{P}(M)$, as required in the hypotheses of Theorem 4.5. A more *ad hoc* study of this particular dual pair using certain transversality properties of the submersions $\pi_1$ and $\pi_2$ shows that the completeness is not actually needed (see [Bl01]) in order to guarantee leaf correspondence. ♦

## 5 Howe pairs and dual pairs

**Definition 5.1** *Let* $(M, \{\cdot, \cdot\})$ *be a Poisson manifold and* $A, B \subset \mathcal{P}_L(M)$ *be two pseudosubgroups of its local Poisson diffeomorphism pseudogroup. We say that the diagram*

$$
\begin{array}{ccc}
 & (M, \{\cdot, \cdot\}) & \\
 & & \\
\pi_A \swarrow & & \searrow \pi_B \\
 & & \\
(M/A, \{\cdot, \cdot\}_{M/A}) & & (M/B, \{\cdot, \cdot\}_{M/B})
\end{array}
$$

*is a **Howe pair** on* $(M, \{\cdot, \cdot\})$ *if the following two conditions are satisfied at the same time:*

$$
\begin{align}
(\pi_A^* C^\infty(M/A))^c &= \pi_B^* C^\infty(M/B) \tag{5.1} \\
(\pi_B^* C^\infty(M/B))^c &= \pi_A^* C^\infty(M/A) \tag{5.2}
\end{align}
$$

*The superscript* $c$ *in the previous equalities means centralizer with respect to the algebra structure on* $C^\infty(M)$ *given by the Poisson bracket on* $M$. *As in the case of the dual pairs, the Poisson manifold* $(M, \{\cdot, \cdot\})$ *will be called the **equivalence bimodule** of the Howe pair.*

**Proposition 5.2** *Let* $(M, \{\cdot, \cdot\})$ *be a smooth Poisson manifold,* $A \subset \mathcal{P}_L(M)$, *and* $A'$ *its dual. Then, if* $A$ *has the extension property and* $A'$ *is integrable we have that*

$$
\begin{align}
(\pi_A^* C^\infty(M/A))^c &= \pi_{A'}^* C^\infty(M/A') \tag{5.3} \\
(\pi_{A'}^* C^\infty(M/A'))^c &\supset \pi_A^* C^\infty(M/A) \tag{5.4}
\end{align}
$$

*Moreover, if* $A$ *is von Neumann and* $A'$ *has the extension property then the diagram* $(M/A, \{\cdot, \cdot\}_{M/A}) \overset{\pi_A}{\longleftarrow} (M, \{\cdot, \cdot\}) \overset{\pi_{A'}}{\longrightarrow} (M/A', \{\cdot, \cdot\}_{M/A'})$ *is a Howe pair.*

**Proof.** We first establish (5.3), which is equivalent to proving that $\left(C^\infty(M)^A\right)^c = C^\infty(M)^{A'}$: let $f \in C^\infty(M)$ arbitrary, and $g \in C^\infty(M)^A$ an $A$-invariant function with associated Hamiltonian flow $F_t$. Then, for any $m \in M$

$$
\frac{d}{dt} f(F_t(m)) = \mathbf{d}f(F_t(m)) \cdot X_g(F_t(m)) = \{f, g\}(F_t(m)). \tag{5.5}
$$

Now, if $f \in \left(C^\infty(M)^A\right)^c$ then $\{f, g\} = 0$ in (5.5) and therefore $f \circ F_t(m) = f(m)$. Since the $A$–invariant function $g$ and the point $m$ are arbitrary, and $A$ has the extension property, we can conclude that $f \in C^\infty(M)^{A'}$. Conversely, if $f \in C^\infty(M)^{A'}$, then $f \circ F_t = f$ and therefore (5.5) implies that $f \in \left(C^\infty(M)^A\right)^c$. Expression (5.4) can be obtained by taking centralizers on both sides of (5.3).

Suppose now that $A$ is von Neumann. In order to conclude that we have a Howe pair we just need to show that $(\pi_{A'}^* C^\infty(M/A'))^c \subset \pi_A^* C^\infty(M/A)$ or, equivalently, that $\left(C^\infty(M)^{A'}\right)^c \subset C^\infty(M)^A$. Let $f \in \left(C^\infty(M)^{A'}\right)^c$. Since $A'$ has the extension property, any element $\mathcal{F}_T \in G_{A''}$ can be written as the finite composition of locally defined flows $F_t$ associated to the Hamiltonian vector fields of $A'$-invariant globally defined functions $h \in C^\infty(M)^{A'}$. Given that for any of those functions we have that $\{f, h\} = 0$, it is clear that $f \circ F_t = f|_{\mathrm{Dom}(F_t)}$. Now, the von Neumann character of $A$ implies that for any $\phi \in A$ and $m \in M$ arbitrary, there exists $\mathcal{F}_T \in G_{A''}$ such that $f \circ \phi(m) = f \circ \mathcal{F}_T(m) = f(m)$, which guarantees the $A$-invariance of $f$.    ∎

**Corollary 5.3** *Let* $(M, \{\cdot, \cdot\})$ *be a smooth Poisson manifold and* $A, B \in \mathcal{P}_L(M)$ *be two pseudosubgroups of its local Poisson diffeomorphism pseudogroup* $\mathcal{P}_L(M)$ *that have the extension property. If the diagram* $(M/A, \{\cdot, \cdot\}_{M/A}) \xleftarrow{\pi_A} (M, \omega) \xrightarrow{\pi_B} (M/B, \{\cdot, \cdot\}_{M/B})$ *is a dual pair then it is also a Howe pair.*

As a corollary to the previous result we can easily obtain the following well known fact:

**Corollary 5.4** *If the diagram* $(P_1, \{\cdot, \cdot\}_{P_1}) \xleftarrow{\pi_1} (M, \omega) \xrightarrow{\pi_2} (P_2, \{\cdot, \cdot\}_{P_2})$ *is a Lie–Weinstein dual pair with connected fibers then it is a Howe pair.*

**Proof.**    It is a straightforward consequence of Example 3.7 where we saw that any Lie–Weinstein dual pair with connected fibers can be understood as a dual pair in our sense with respect to two pseudosubgroups $G_\mathcal{A}, G_\mathcal{B} \subset \mathcal{P}_L(M)$. These pseudosubgroups have the extension property by Lemma 3.6 and hence the hypotheses of the previous corollary are satisfied in our case.    ∎

Even though Corollary 5.3 shows that in the presence of the extension property any dual pair is a Howe pair, the following example demonstrates that the converse is not, in general, true.

**Example 5.5 A Howe pair that is not a dual pair**: Let $(\mathbb{T}^2, \omega)$ be the two torus thought of as a symplectic manifold with the form $\omega$ given by the standard area form. Consider a Poisson action of the additive group $(\mathbb{R}, +)$ on $\mathbb{T}^2$ via an irrational flow. It is straightforward to check that the Poisson diffeomorphisms group $A_\mathbb{R} \subset \mathcal{P}(\mathbb{T}^2)$ associated to this action generates a Howe pair $\mathbb{T}^2/A_\mathbb{R} \xleftarrow{\pi_{A_\mathbb{R}}} \mathbb{T}^2 \xrightarrow{\pi_{A_\mathbb{R}'}} \mathbb{T}^2/A_\mathbb{R}'$ that is not a dual pair. Indeed, notice first that the only $A_\mathbb{R}$–invariant open subset of $\mathbb{T}^2$ is $\mathbb{T}^2$ itself hence, since $C^\infty(\mathbb{T}^2)^{A_\mathbb{R}}$ is made of constant functions the dual distribution $A_\mathbb{R}'$ is trivial and $C^\infty(\mathbb{T}^2)^{A_\mathbb{R}'} = C^\infty(\mathbb{T}^2)$, necessarily. It is clear that in these circumstances $(C^\infty(\mathbb{T}^2)^{A_\mathbb{R}})^c = C^\infty(\mathbb{T}^2)^{A_\mathbb{R}'}$ and $(C^\infty(\mathbb{T}^2)^{A_\mathbb{R}'})^c = C^\infty(\mathbb{T}^2)^{A_\mathbb{R}}$. Nevertheless, the orbits of the $A_\mathbb{R}$–action are strictly contained inside the only leaf of the distribution $A_\mathbb{R}''$, which implies that $A_\mathbb{R}$ is not von Neumann and thereby does not generate a dual pair.

This example also shows that Howe's condition is not enough to ensure symplectic leaf correspondence. Indeed, the remarks in the preceding paragraph indicate that the Howe pair associated to $A_\mathbb{R}$ is $\mathbb{T}^2/A_\mathbb{R} \xleftarrow{\pi_{A_\mathbb{R}}} \mathbb{T}^2 \xrightarrow{\mathrm{id}} \mathbb{T}^2$. Now, the right hand side leg of this pair has just one symplectic leaf (the entire two torus $\mathbb{T}^2$) while, for the left hand side, every point in $\mathbb{T}^2/A_\mathbb{R}$ is a symplectic leaf since $C^\infty(\mathbb{T}^2/A_\mathbb{R})$ consists of constant functions.    ♦

# 6   Hamiltonian Poisson subgroups

In this section we will study the properties of the diagrams $(M/A, \{\cdot, \cdot\}_{M/A}) \xleftarrow{\pi_A} (M, \{\cdot, \cdot\}) \xrightarrow{\pi_{A'}} (M/A', \{\cdot, \cdot\}_{M/A'})$ induced by weakly and strongly Hamiltonian subgroups $A \subset \mathcal{P}(M)$. Since we are dealing with actual subgroups of $\mathcal{P}(M)$, Proposition 3.8 guarantees the integrability of the polar distribution $A'$ which we will not need to put as a hypothesis.

In Example 5.5 we identified a weakly Hamiltonian subgroup that induced a Howe pair. In our first result in this section, Proposition 6.1, we will show that this is not a coincidence since any weakly Hamiltonian subgroup endowed with the extension property always has a Howe pair associated. We also saw in that example that the (weak) Hamiltonian condition is not sufficient to generate a dual pair; in Proposition 6.2 we will show that if we add to the Hamiltonian hypothesis the property of separation $A$–orbits then we are guaranteed to obtain a dual pair.

**Proposition 6.1** *Let $(M, \{\cdot, \cdot\})$ be a Poisson manifold and $A \subset \mathcal{P}(M)$ be a weakly Hamiltonian subgroup of its Poisson diffeomorphism group. If $A$ has the extension property then the diagram $(M/A, \{\cdot, \cdot\}_{M/A}) \xleftarrow{\pi_A} (M, \{\cdot, \cdot\}) \xrightarrow{\pi_{A'}} (M/A', \{\cdot, \cdot\}_{M/A'})$ is a Howe pair.*

**Proof.** By Proposition 3.8 the polar distribution $A'$ is always integrable in this case. The conclusions of Proposition 5.2 show that we just need to prove that $(\pi_{A'}^* C^\infty(M/A'))^c \subset \pi_A^* C^\infty(M/A)$. Hence, let $\phi \in A$ and $m \in M$ arbitrary. Since, by hypothesis, the group $A$ is weakly Hamiltonian, $\phi(m)$ can be written as $\phi(m) = F_{t_1}^1 \circ F_{t_2}^2 \circ \cdots \circ F_{t_k}^k(m)$, with $F_{t_i}^i$ the flow of a Hamiltonian vector field $X_{h_i}$ associated to a function $h_i$ in the centralizer $\left(C^\infty(M)^A\right)^c$. We assume for the sake of simplicity that $\phi(m) = F_t(m)$, with $F_t$ the flow of $X_h$, $h \in \left(C^\infty(M)^A\right)^c$. Due to the expression (5.3), the function $h$ can be written as $h = g \circ \pi_{A'}$, with $g \in C^\infty(M/A')$. Let now $f \in (\pi_{A'}^* C^\infty(M/A'))^c$. Given that $\{f, h\} = \{f, g \circ \pi_{A'}\} = 0$ we can conclude that $f \circ \phi(m) = f \circ F_t(m) = f(m)$. As we can reproduce this process for any $\phi \in A$ and $m \in M$ we have that $f \in C^\infty(M)^A = \pi_A^* C^\infty(M/A)$, as required. ∎

**Proposition 6.2** *Let $(M, \{\cdot, \cdot\})$ be a Poisson manifold and $A \subset \mathcal{P}(M)$ be a subgroup of its Poisson diffeomorphism group.*

**(i)** *If $A$ is strongly (resp. weakly) Hamiltonian and has the extension property, then*

$$A \subset G_{\mathcal{A}''} \qquad (resp. \ A \cdot m \subset G_{\mathcal{A}''} \cdot m \ for \ any \ m \in M).$$

**(ii)** *If $C^\infty(M)^A$ separates the $A$–orbits on $M$ then, for any $m \in M$, we have that*

$$G_{\mathcal{A}''} \cdot m \subset A \cdot m.$$

**(iii)** *If $A$ is (strongly or weakly) Hamiltonian and has the extension property, and $C^\infty(M)^A$ separates the $A$–orbits on $M$, then $A$ is von Neumann and the diagram $(M/A, \{\cdot, \cdot\}_{M/A}) \xleftarrow{\pi_A} (M, \{\cdot, \cdot\}) \xrightarrow{\pi_{A'}} (M/A', \{\cdot, \cdot\}_{M/A'})$ is a dual pair. Additionally, if $A'$ has the extension property it is also a Howe pair.*

**Proof.** **(i)** Let $\phi \in A$ be arbitrary. Since $A$ is strongly (resp. weakly) Hamiltonian, $\phi$ (resp. $\phi(m)$ for any $m \in M$) can be written as $\phi = F_{t_1}^1 \circ F_{t_2}^2 \circ \cdots \circ F_{t_k}^k$ (resp. $\phi(m) = F_{t_1}^1 \circ F_{t_2}^2 \circ \cdots \circ F_{t_k}^k(m)$), with $F_{t_i}^i$ the flow of a Hamiltonian vector field $X_{h_i}$ associated to a function $h_i$ in the centralizer $\left(C^\infty(M)^A\right)^c$. In order to keep the exposition simple we assume that $\phi = F_t$, with $F_t$ the flow of $X_h$, $h \in \left(C^\infty(M)^A\right)^c$. Due to (5.3), the function $h$ can be written as $h = l \circ \pi_{A'}$, with $l \in C^\infty(M/A')$. Consequently, $X_h = X_{l \circ \pi_{A'}}$ and hence $F_t = \phi \in G_{\mathcal{A}''}$ (resp. $F_t(m) = \phi(m) \in G_{\mathcal{A}''} \cdot m$), as required.

**(ii)** Any element in $\mathcal{F}_T \in G_{\mathcal{A}''}$ can be written as a finite composition of Hamiltonian flows $F_t$ associated to functions $f \circ \pi_{A'}|_U$, $f \in C^\infty(U/A')$, $U$ an open $A'$–invariant set. Then, for any $h \in C^\infty(M)^A$ and any $m \in U$ we have that

$$\frac{d}{dt} h\left(F_t(m)\right) = \{h|_U, f \circ \pi_{A'}|_U\}_U \left(F_t(m)\right) = -\mathbf{d}(f \circ \pi_{A'}|_U)\left(F_t(m)\right) \cdot X_h\left(F_t(m)\right) = 0,$$

that is, any function $h \in C^\infty(M)^A$ is constant along the Hamiltonian flow of $f \circ \pi_{A'}|_U$. Now, since $C^\infty(M)^A$ separates the $A$–orbits on $M$, we can conclude that, for any point $m \in M$, the set $F_t(m)$ is included in a single $A$-orbit, namely, $F_t(m) \subset A \cdot m$ and therefore $G_{\mathcal{A}''} \cdot m \subset A \cdot m$, as required.

**(iii)** Parts **(i)** and **(ii)** imply in the context of our hypotheses that for any $m \in M$, $A \cdot m = G_{\mathcal{A}''} \cdot m''$ and, consequently, $M/A = M/A''$. This proves that $A$ is von Neumann and therefore that the diagram $(M/A, \{\cdot, \cdot\}_{M/A}) \xleftarrow{\pi_A} (M, \{\cdot, \cdot\}) \xrightarrow{\pi_{A'}} (M/A', \{\cdot, \cdot\}_{M/A'})$ is a dual pair. Corollary 5.3 ensures that it is also a Howe pair in the presence of the extension property for $A'$. ∎

# 7 Dual pairs induced by canonical Lie group actions

In this section we will analyze under what circumstances we can construct von Neumann and Howe pairs using the subgroups $A_G := \{\Phi_g \mid g \in G\}$ of the Poisson diffeomorphism group $\mathcal{P}(M)$ associated to the canonical action $\Phi : G \times M \to M$ of a Lie group $G$ on a Poisson manifold $(M, \{\cdot, \cdot\})$. Recall that in this setup, as we already mentioned in Example 3.13, the polar distribution $A'_G$ is always integrable (Proposition 3.8) and the projection onto the corresponding leaf space $\mathcal{J} : M \to M/A'_G$ is referred to as the ***optimal momentum map***. The reason behind this denomination is (see [OR02] for the details) that the Hamiltonian flow $F_t$ associated to any $G$–invariant function $f \in C^\infty(U)^G$ defined on any $G$–invariant open subset $U$ of $M$ preserves the level sets of $\mathcal{J}$, that is, $\mathcal{J} \circ F_t = \mathcal{J}$ (Noether's Theorem). Moreover, by construction, the level sets of this map are the smallest submanifolds of $M$ preserved by $G$–equivariant Hamiltonian flows on $M$. Also, the map $\mathcal{J}$ is universal in the category of the momentum maps that can be associated to the $G$–symmetry of $(M, \{\cdot, \cdot\})$.

## 7.1 Properness, Hamiltonian actions, and dual pairs

An action $\Phi : G \times M \to M$ of a Lie group $G$ on a manifold $M$ is said to be ***proper*** if the following condition is satisfied: given two convergent sequences, $\{m_n\}$ and $\{g_n \cdot m_n\}$ in $M$, there exists a convergent subsequence $\{g_{n_k}\}$ in $G$. The following proposition is a summary of well known facts about proper actions that we will use in the sequel. The reader is encouraged to check with [ACG91, DuKo99, OR02] for proofs.

**Proposition 7.1** *Let $\Phi : G \times M \to M$ be a proper action of a Lie group $G$ on a smooth manifold $M$. Then:*

**(i)** $A_G := \{\Phi_g \mid g \in G\}$ *has the extension property.*

**(ii)** $C^\infty(M)^{A_G} = C^\infty(M)^G$ *separates the $G$–orbits.*

**(iii)** *The isotropy subgroup $G_m$ of any point $m \in M$ is compact.*

**Theorem 7.2** *Let $G$ be a Lie group acting canonically and properly on the Poisson manifold $(M, \{\cdot, \cdot\})$ via the map $\Phi : G \times M \to M$. Let $A_G \subset \mathcal{P}(M)$ be the subgroup of $\mathcal{P}(M)$ defined by $A_G := \{\Phi_g : M \to M \mid g \in G\}$ and $A'_G$ its polar. Let $\pi : M \to M/A_G$ be the canonical projection of $M$ onto the quotient $M/A_G$ and $\mathcal{J} : M \to M/A'_G$ be the associated optimal momentum map. If $A_G$ is (strongly or weakly) Hamiltonian then it is von Neumann and therefore the diagram $(M/A_G, \{\cdot, \cdot\}_{M/A_G}) \xleftarrow{\pi} (M, \{\cdot, \cdot\}) \xrightarrow{\mathcal{J}} (M/A'_G, \{\cdot, \cdot\}_{M/A'_G})$ is a dual pair.*

**Proof.** It is a straightforward consequence of Propositions 6.2 and 7.1. ∎

**Corollary 7.3** *In the same setup as in the previous theorem, if $A_G$ is (strongly or weakly) Hamiltonian then, the $G$–orbits are connected and path connected.*

**Proof.** The condition on $A_G$ being Hamiltonian implies via the previous theorem that $A_G$ is von Neumann and therefore, for any $m \in M$ the orbit $G \cdot m$ equals $G_{\mathcal{A}''} \cdot m$ which is connected and path connected. ∎

Theorem 7.2 shows that properness in a canonical $G$–action is a condition that added to the Hamiltonian character is sufficient to ensure that the corresponding transformation group $A_G \subset \mathcal{P}(M)$ is von Neumann. However, as the following example shows, this condition is not necessary.

**Example 7.4 The coadjoint action produces von Neumann subgroups of $\mathcal{P}(\mathfrak{g}^*)$.** Let $G$ be a connected Lie group, $\mathfrak{g}$ its Lie algebra, and $\mathfrak{g}^*$ its dual. Let $\{\cdot, \cdot\}$ be the $+$–Lie-Poisson bracket that makes $\mathfrak{g}^*$ into a Poisson manifold. More specifically, for any $f, h \in C^\infty(\mathfrak{g}^*)$ and $\mu \in \mathfrak{g}^*$, we define $\{f, h\}(\mu) := \left\langle \mu, \left[ \frac{\delta f}{\delta \mu}, \frac{\delta h}{\delta \mu} \right] \right\rangle$, where the symbol $\langle \cdot, \cdot \rangle$ denotes the natural pairing of $\mathfrak{g}^*$ with $\mathfrak{g}$ and the elements $\frac{\delta f}{\delta \mu}, \frac{\delta h}{\delta \mu} \in \mathfrak{g}$ are determined by the expressions $Df(\mu) \cdot \rho := \left\langle \rho, \frac{\delta f}{\delta \mu} \right\rangle, Dh(\mu) \cdot \rho := \left\langle \rho, \frac{\delta h}{\delta \mu} \right\rangle$, for all $\rho \in \mathfrak{g}^*$. Given that for any $f \in C^\infty(\mathfrak{g}^*)$, $g \in G$, and $\mu \in \mathfrak{g}^*$ we have that

$$\frac{\delta f}{\delta \left( \mathrm{Ad}^*_{g^{-1}} \mu \right)} = \mathrm{Ad}_g \left( \frac{\delta \left( f \circ \mathrm{Ad}^*_{g^{-1}} \right)}{\delta \mu} \right),$$

it can be readily verified that coadjoint action of $G$ on $\mathfrak{g}^*$ is canonical and has the identity as standard momentum map associated. Also, for any $f \in C^\infty(\mathfrak{g}^*)$ and $\mu \in \mathfrak{g}^*$ $X_f(\mu) = -\mathrm{ad}^*_{\frac{\delta f}{\delta \mu}} \mu$.

We now check that $A_G$ is von Neumann. Let $\mu \in \mathfrak{g}^*$ be arbitrary and $U \subset \mathfrak{g}^*$ be an open $G$–invariant neighborhood of the coadjoint orbit of the element $\mu$. Let $f \in C^\infty(U)^G$. Then for any $\xi \in \mathfrak{g}$ and $\rho \in U$ we have that

$$\langle X_f(\rho), \xi \rangle = -\left\langle \mathrm{ad}^*_{\frac{\delta f}{\delta \rho}} \rho, \xi \right\rangle = -\left\langle \rho, \left[ \frac{\delta f}{\delta \rho}, \xi \right] \right\rangle = \frac{d}{dt} \bigg|_{t=0} \left\langle \rho, \mathrm{Ad}_{\exp t\xi} \frac{\delta f}{\delta \rho} \right\rangle$$
$$= \frac{d}{dt} \bigg|_{t=0} \left\langle \mathrm{Ad}^*_{\exp t\xi} \rho, \frac{\delta f}{\delta \rho} \right\rangle = \left\langle \mathrm{ad}^*_\xi \rho, \frac{\delta f}{\delta \rho} \right\rangle = \frac{d}{dt} \bigg|_{t=0} f \left( \mathrm{Ad}^*_{\exp t\xi} \rho \right) = 0,$$

where the last equality follows from the $G$–invariance of the function $f$. This computation shows that $A'_G(\mu) = \{0\}$ for all $\mu \in \mathfrak{g}^*$. The connectedness of the group $G$ automatically implies that $A''_G = A_G$ and therefore $A_G$ is von Neumann.

The symplectic leaf correspondence for the legs of the diagram $\mathfrak{g}^*/G \leftarrow \mathfrak{g}^* \to \mathfrak{g}^*$ guaranteed in this case by Theorem 4.5 is a restatement of the well known fact that the symplectic leaves of $(\mathfrak{g}^*, \{\cdot, \cdot\})$ are the coadjoint orbits. ♦

## 7.2 Compact connected Lie group symplectic actions and von Neumann subgroups of $\mathcal{P}(M)$

**Definition 7.5** *Let $G$ be a compact connected Lie group with Lie algebra $\mathfrak{g}$ acting canonically on the symplectic manifold $(M, \omega)$ via the map $\Phi : G \times M \to M$. Let $A_G \subset \mathcal{P}(M)$ be the subgroup of $\mathcal{P}(M)$ defined by $A_G := \{\Phi_g : M \to M \mid g \in G\}$. Let $\xi \in \mathfrak{g}$ and $T(\xi)$ be the torus defined by*

$$T(\xi) := \overline{\{\exp t\xi \mid t \in \mathbb{R}\}}.$$

We will say that the element $\xi$ has a **coisotropic torus** associated when the orbits of the $T(\xi)$–action on $M$ are coisotropic.

**Theorem 7.6** *Let $G$ be a compact connected Lie group with Lie algebra $\mathfrak{g}$ acting canonically on the symplectic manifold $(M, \omega)$ via the map $\Phi : G \times M \to M$. Let $A_G \subset \mathcal{P}(M)$ be the subgroup of $\mathcal{P}(M)$ defined by $A_G := \{\Phi_g : M \to M \mid g \in G\}$ and $A'_G$ its polar. Let $\pi : M \to M/A_G$ be the canonical projection of $M$ onto the quotient $M/A_G$ and $\mathfrak{J} : M \to M/A'_G$ be the associated optimal momentum map. Let $\mathbb{T}$ be a maximal torus of $G$ and suppose that its Lie algebra $\mathfrak{t}$ has a basis $\{\xi_1, \ldots, \xi_k\}$ whose elements have coisotropic tori $T(\xi_i)$ associated. Then, $A_G$ is weakly Hamiltonian and von Neumann.*

**Proof.**    Since the action of any compact group is always proper, according to Theorem 7.2, it suffices to prove that $A_G$ is weakly Hamiltonian, which will be a consequence of the following lemma:

**Lemma 7.7** *Suppose that we have the same setup as Theorem 7.6. Then, for any $\xi \in \mathfrak{g}$ that has a coisotropic torus $T(\xi)$ associated and any $m \in M$, there exists a function $f \in (C^\infty(M)^{A_G})^c$ such that if $F_t$ is the flow of the Hamiltonian vector field $X_f$ then:*

$$\exp \xi \cdot m = F_1(m).$$

**Proof.**    Let $\xi_M \in \mathfrak{X}(M)$ be the vector field that assigns to any point $m \in M$, the infinitesimal generator at $m$ associated to the element $\xi \in \mathfrak{g}$. The canonical character of the action implies that

$$0 = \mathcal{L}_{\xi_M} \omega = \mathbf{i}_{\xi_M} \mathbf{d}\omega + \mathbf{d}\left(\mathbf{i}_{\xi_M} \omega\right) = \mathbf{d}\left(\mathbf{i}_{\xi_M} \omega\right),$$

that is, the one form $\alpha := \mathbf{i}_{\xi_M} \omega$ is closed. Consider now the subsets of $G$ defined by $K := \{\exp t\xi \mid t \in \mathbb{R}\}$ and $T(\xi) := \bar{K}$, where the bar over $K$ means closure. As we already pointed out subset $T(\xi)$ is a closed connected Abelian subgroup of $G$ and therefore a torus. Notice that for any $m \in M$ we have that $T(\xi) \cdot m \subset \overline{K \cdot m}$; indeed, if $t \cdot m \in T(\xi) \cdot m$, there exists a sequence $\{k_n\} \subset K$ of elements in $K$ such that $k_n \to t$, which implies that $k_n \cdot m \to t \cdot m$ and therefore $t \cdot m \in \overline{K \cdot m}$. Hence, since the restriction $\alpha|_{K \cdot m} = 0$ we have that $\alpha|_{\overline{K \cdot m}} = 0$, and therefore $\alpha|_{T(\xi) \cdot m} = 0$. By the Relative Poincaré Lemma (see for instance Corollary 7.5 in page 362 of [LM87]) there exists a neighborhood $U$ of $T(\xi) \cdot m$, which by the compactness of $T(\xi)$ can be chosen $T(\xi)$–invariant, and a function $h \in C^\infty(U)$ such that $\mathbf{d}h|_U = \alpha|_U$. This statement amounts to saying that the function $h \in C^\infty(U)$ is a momentum map for the canonical action of $K$ on the symplectic manifold $(U, \omega|_U)$.

Now, by shrinking $U$ if necessary and using the hypothesis on the coisotropic character of the torus $T(\xi)$, we can represent $U$ by a normal form coordinate chart around the point $m$ similar to the ones introduced in the Appendix (Theorem 8.1), that is, we can assume without loss of generality that

$$U \cong T(\xi) \times_{T(\xi)_m} \mathfrak{m}_r^*,$$

where $\mathfrak{m}$ is a $\mathrm{Ad}_{T(\xi)_m}$-invariant complement to the Lie algebra $\mathrm{Lie}\left(T(\xi)_m\right)$ in $\mathfrak{k} := \{\eta \in \mathrm{Lie}\left(T(\xi)\right) \mid \eta_M(m) \in \left(\mathrm{Lie}\left(T(\xi)\right) \cdot m\right)^\omega\}$. The point $m$ is represented in these coordinates by $[e, 0]$, and $\mathfrak{m}_r^*$ is a $T(\xi)_m$-equivariant ball of radius $r > 0$ small enough centered at the origin of $\mathfrak{m}^*$. Let $\phi_r : \mathfrak{m}^* \to \mathbb{R}$ be a smooth, $T(\xi)_m$–invariant, and compactly supported function such that $\phi_r(\eta) = 0$, for any $\eta \in \mathfrak{m}^* \setminus \mathfrak{m}_r^*$, and $\phi_r(W) = 1$ for a $T(\xi)_m$-invariant neighborhood $W \subset \mathfrak{m}_r^*$. Let $\Phi$ be the $T(\xi)$–invariant function defined by

$$\Phi : \quad U \cong T(\xi) \times_{T(\xi)_m} \mathfrak{m}_r^* \quad \longrightarrow \quad \mathbb{R}$$
$$[k, \eta] \quad \longmapsto \quad \phi_r(\eta).$$

Notice that $\Phi$ is zero off the open $T(\xi)$–invariant set $U$ and therefore it can be trivially extended to a $T(\xi)$–invariant function, we will call it equally $\Phi \in C^\infty(M)^{T(\xi)}$, on the entire space. The reconstruction equations (8.5) through (8.7) applied to $\Phi$ (use the Abelian character of $T(\xi)$) imply that the Hamiltonian vector field $X_\Phi$ equals

$$X_\Phi(z) = \begin{cases} (D_{\mathbf{m}^*}\phi_r)_M(m) & \text{if} \quad m \in U \\ 0 & \text{if} \quad m \in M \setminus U. \end{cases}$$

Let $f = \Phi h$. Now, given that $X_f = X_{\Phi h} = \Phi X_h + h X_\Phi$ and $\Phi$ is constant in the $T(\xi)$–invariant neighborhood $N \simeq T(\xi) \times_{T(\xi)_m} W$ around $m$, we have that $X_f(z) = X_h(z) = \xi_M(z)$ for any $z \in N$. Consequently, if $F_t$ is the flow of the vector field $X_f$, it is clear that

$$\exp\xi \cdot m = F_1(m).$$

In order to finish the proof we just need to show that $f \in (C^\infty(M)^{A_G})^c$. This is indeed so because for any $G$-invariant function $l \in C^\infty(M)^{A_G}$ we have that $\{f, l\}(z) = \mathbf{d}f(z) \cdot X_l(z) = 0$ for any $z \in M \setminus U$. Also, when $z \in U$

$$\begin{aligned} \{f, l\}(z) &= \{\Phi h, l\}(z) = \Phi(z)\{h, l\}(z) + h(z)\{\Phi, l\}(z) \\ &= -\Phi(z)\left(\mathbf{d}l(z) \cdot X_h(z)\right) - h(z)\left(\mathbf{d}l(z) \cdot X_\Phi(z)\right) \\ &= -\Phi(z)\left(\mathbf{d}l(z) \cdot \xi_M(z)\right) - h(z)\left(\mathbf{d}l(z) \cdot (D_{\mathbf{m}^*}\phi_r)_M(z)\right) = 0, \end{aligned}$$

due to the $G$-invariance of $l$. Hence $\{f, l\} = 0$ for any $l \in C^\infty(M)^{A_G}$ and, consequently, $f \in (\pi^* C^\infty(M/A_G))^c$, as required.  ▼

We conclude the proof of the theorem by noting that since the group $G$ is compact and connected, any element $g \in G$ can be written as $g = hlh^{-1}$, with $l \in \mathbb{T}$. As $\mathbb{T}$ is Abelian and connected, there exist real numbers $t_1, \ldots, t_k$ such that $l = \exp t_1\xi_1 \cdots \exp t_k\xi_k$. Hence, for any $m \in M$ we can write

$$\begin{aligned} g \cdot m &= h\exp t_1\xi_1 \cdots \exp t_k\xi_k h^{-1} \cdot m \\ &= h\exp t_1\xi_1 h^{-1} \cdots h\exp t_k\xi_k h^{-1} \cdot m = \exp t_1(\mathrm{Ad}_h\xi_1) \cdots \exp t_k(\mathrm{Ad}_h\xi_k) \cdot m \end{aligned}$$

A straightforward computation shows that $T(\mathrm{Ad}_h(\xi_i)) = hT(\xi_i)h^{-1}$ and that, as a consequence, the coisotropy of the torus $T(\xi_i)$ implies that of $T(\mathrm{Ad}_h(\xi_i))$. Therefore, by the previous lemma we have that

$$g \cdot m = F_1^1 \circ \cdots \circ F_1^k(m),$$

with each $F_t^i$ the flow of a Hamiltonian vector field $X_{f_i}$ associated to a function $f_i \in (\pi^* C^\infty(M/A_G))^c$.  ∎

The reader may be wondering if the coisotropy hypothesis in the statement of Theorem 7.6 is not just a technical requirement that appears in the proof of Lemma 7.7 and that could be eliminated by using different techniques in the proof. The following example, that I owe to J. Montaldi and T. Tokieda, shows that this is not the case, that is, in the absence of additional hypotheses, compactness and connectedness in the Lie group $G$ associated to a symplectic action do not suffice to ensure that the corresponding transformation group $A_G$ is von Neumann.

**Example 7.8 A compact and connected canonical group action that is not von Neumann.** Let $M := \mathbb{T}^2 \times \mathbb{T}^2$ be the product of two tori whose elements we will denote by the four–tuples $(e^{i\theta_1}, e^{i\theta_2}, e^{i\psi_1}, e^{i\psi_2})$. We endow $M$ with the symplectic structure $\omega$ defined by $\omega := \mathbf{d}\theta_1 \wedge \mathbf{d}\theta_2 + \sqrt{2}\,\mathbf{d}\psi_1 \wedge \mathbf{d}\psi_2$. We now consider the canonical circle action given by $e^{i\phi} \cdot (e^{i\theta_1}, e^{i\theta_2}, e^{i\psi_1}, e^{i\psi_2}) := (e^{i(\theta_1+\phi)}, e^{i\theta_2}, e^{i(\psi_1+\phi)}, e^{i\psi_2})$. This action does not satisfy the coisotropy hypothesis and, as we will now verify, the associated transformation group $A_{S^1} \subset \mathcal{P}(M)$ is not von Neumann. Indeed, the set $C^\infty(M)^{S^1}$ comprises the functions $f$ of the form $f \equiv f(e^{i\theta_2}, e^{i\psi_2}, e^{i(\theta_1-\psi_1)})$. An inspection of the Hamiltonian flows associated to such functions readily shows that the leaves of $A'_{S^1}$ fill densely the manifold $M$. This implies that $C^\infty(M)^{A'_{S^1}}$ is made up by constant functions and therefore $A''_{S^1}(m) = \{0\}$, for all $m \in M$. Consequently, $A_{S^1}$ is not von Neumann.

Notice that this example shows that the polar of a regular integrable distribution even though it is integrable, it is not, in general, regular. More specifically, even though the projection $\pi_{A_{S^1}} : M \to M/A_{S^1}$ is a surjective submersion, this is not true in the case of $\pi_{A'_{S^1}} : M \to M/A'_{S^1}$. ◆

## 7.3 Tubewise Hamiltonian actions and dual pairs

In the Appendix (see Section 8.3) the reader can find an in depth study of the conditions under which the proper canonical action of a connected Lie group $G$ on a symplectic manifold $(M, \omega)$ is strongly tubewise Hamiltonian. More specifically, in that section it is explained how under some circumstances, for any point $m \in M$ there is an open $G$–invariant neighborhood of its orbit such that the restriction of the $G$–action to this neighborhood has a standard momentum map associated, thus implying that the action is strongly tubewise Hamiltonian. The question that we will try to answer in this section is the following: is there any situation where the strongly tubewise Hamiltonian condition implies that the action is weakly Hamiltonian and therefore induces a dual pair? The following result provides some answers to this question.

**Theorem 7.9** *Let $G$ be a connected Lie group with Lie algebra $\mathfrak{g}$ acting canonically and properly on the symplectic manifold $(M, \omega)$ via the map $\Phi : G \times M \to M$. Let $A_G \subset \mathcal{P}(M)$ be the subgroup of $\mathcal{P}(M)$ defined by $A_G := \{\Phi_g : M \to M \mid g \in G\}$ and $A'_G$ its polar. Let $\pi : M \to M/A_G$ be the canonical projection of $M$ onto the quotient $M/A_G$ and $\mathcal{J} : M \to M/A'_G$ be the associated optimal momentum map. For any $m \in M$ let $\mathfrak{k}_m \subset \mathfrak{g}$ be the Lie subalgebra of $\mathfrak{g}$ defined by*

$$\mathfrak{k}_m = \{\eta \in \mathfrak{g} \mid \eta_M(m) \in (\mathfrak{g} \cdot m)^\omega\},$$

*$K_m \subset G$ be the (unique) connected Lie subgroup generated by it $\mathfrak{k}_m$, and $\gamma_m \in \Omega^1(G; \mathfrak{g}^*)$ be the $G$–equivariant, $\mathfrak{g}^*$–valued one form defined by*

$$\langle \gamma_m(g) \cdot T_e L_g \cdot \eta, \xi \rangle := -\omega(m)\left((\mathrm{Ad}_{g^{-1}}\xi)_M(m), \eta_M(m)\right) \quad \text{for any} \quad g \in G, \, \xi, \eta \in \mathfrak{g}.$$

*Suppose that for any $m \in M$, the orbit $G \cdot m$ is coisotropic, there exists a $\mathrm{Ad}(K_m)$-invariant complement to $\mathfrak{k}_m$ in $\mathfrak{g}$, $\mathfrak{k}_m$ is Abelian, and $\gamma_m$ is exact (which happens for instance when $H^1(G) = 0$ or when the orbit $G \cdot m$ is isotropic). Then, $A_G$ is weakly Hamiltonian and von Neumann.*

**Proof.** We will show that in the presence of our hypotheses a conclusion similar to that of Lemma 7.7 holds, that is, we will see that for any $\xi \in \mathfrak{g}$ and any $m \in M$, there exists a function $f \in (C^\infty(M)^{A_G})^c$ such that if $F_t$ is the flow of the Hamiltonian vector field $X_f$ then,

$$\exp \xi \cdot m = F_1(m).$$

Indeed, for a fixed $m \in M$, the exactness of $\gamma_m$ guarantees, by Proposition 8.2, that there exists a $G$–invariant neighborhood $U$ of the orbit $G \cdot m$ where the restriction of the $G$–action has a standard momentum map associated and, consequently, for any $\xi \in \mathfrak{g}$ we have that $\exp \xi \cdot m = F_1(m)$, with $F_t$ the flow of the Hamiltonian vector field in $U$ associated to a function $h \in C^\infty(U)$ that can be constructed by taking the $\xi$–component of the tubular momentum map.

We now proceed in a way that mimics the proof of Lemma 7.7. First, by shrinking $U$ if necessary, we can represent it by a normal form coordinate chart around the point $m$, that is, we can assume without loss of generality that $U \cong G \times_{G_m} \mathfrak{m}_r^*$, where $\mathfrak{m}$ is a $\mathrm{Ad}_{G_m}$-invariant complement to the Lie algebra $\mathrm{Lie}(G_m)$ in $\mathfrak{m}_m$. The point $m$ is represented in these coordinates by $[e, 0]$, and $\mathfrak{m}_r^*$ is an open $G_m$–equivariant ball of radius $r > 0$ small enough centered at the origin of $\mathfrak{m}^*$. Let $\phi_r : \mathfrak{m}^* \to \mathbb{R}$ be a smooth $G_m$–invariant compactly supported function such that $\phi_r(\eta) = 0$, for any $\eta \in \mathfrak{m}^* \setminus \mathfrak{m}_r^*$, and $\phi_r(W) = 1$ for a $G_m$–invariant neighborhood $W \subset \mathfrak{m}_r^*$. Let $\Phi$ be the $G$-invariant function defined by

$$\Phi : \quad U \cong G \times_{G_m} \mathfrak{m}_r^* \quad \longrightarrow \quad \mathbb{R}$$
$$[k, \eta] \quad \longmapsto \quad \phi_r(\eta).$$

Notice that $\Phi$ is zero off the open $G$–invariant set $U$ and therefore can be trivially extended to a $G$–invariant function, we will call it equally $\Phi \in C^\infty(M)^G$, on the entire space. The reconstruction equations (8.5) through (8.7) applied to $\Phi$ imply that the Hamiltonian vector field $X_\Phi$ equals

$$X_\Phi(z) = \begin{cases} (D_{\mathfrak{m}^*}\phi_r)_M(m) & \text{if} \quad m \in U \\ 0 & \text{if} \quad m \in M \setminus U. \end{cases} \tag{7.1}$$

Indeed, the hypothesis on the existence of a $\mathrm{Ad}(K_m)$-invariant complement to $\mathfrak{k}_m$ in $\mathfrak{g}$ implies that the map $F$ in (8.4) reduces to $F(\xi, \lambda, \tau) = \mathbb{P}_{\mathfrak{q}^*}(\mathrm{ad}_\tau^* \lambda) + \langle \tau, \cdot \rangle_{\mathfrak{q}}$, whose unique solution for $\tau$ is $\tau \equiv 0$. This implies that the map $\psi \equiv 0$ and, given that by hypothesis the Lie algebra $\mathfrak{k}_m$ is Abelian, the reconstruction equation (8.6) vanishes, thus justifying (7.1).

Let $f = \Phi h$. Now, given that $X_f = X_{\Phi h} = \Phi X_h + h X_\Phi$ and $\Phi$ is constant on the $G$-invariant neighborhood $N \simeq G \times_{G_m} W$ around $m$, we have that $X_f(z) = X_h(z) = \xi_M(z)$ for any $z \in N$. Consequently, if $F_t$ is the flow of the vector field $X_f$, it is clear that

$$\exp \xi \cdot m = F_1(m).$$

In order to finish the proof we just need to show that $f \in (C^\infty(M)^{A_G})^c$. This is indeed so because for any $G$–invariant function $l \in C^\infty(M)^{A_G}$ we have that $\{f, l\}(z) = \mathbf{d}f(z) \cdot X_l(z) = 0$ for any $z \in M \setminus U$. Also, when $z \in U$

$$\begin{aligned}
\{f, l\}(z) &= \{\Phi h, l\}(z) = \Phi(z)\{h, l\}(z) + h(z)\{\Phi, l\}(z) \\
&= -\Phi(z)\left(\mathbf{d}l(z) \cdot X_h(z)\right) - h(z)\left(\mathbf{d}l(z) \cdot X_\Phi(z)\right) \\
&= -\Phi(z)\left(\mathbf{d}l(z) \cdot \xi_M(z)\right) - h(z)\left(\mathbf{d}l(z) \cdot (D_{\mathfrak{m}^*}\phi)_M(z)\right) = 0,
\end{aligned}$$

due to the $G$-invariance of $l$. Hence $\{f, l\} = 0$ for any $l \in C^\infty(M)^{A_G}$ and, consequently, $f \in (\pi^* C^\infty(M/A_G))^c$, as required. ∎

## 7.4 Complete polar distributions and symplectic leaf correspondence

In this section we will show that the polar distribution $A_G'$ relative to the proper and canonical action of a Lie group $G$ on a symplectic manifold is complete (see Definition 2.4). Therefore the quotient space $M/A_G'$ can be written as the orbit space $M/G_{A_c'}$ relative to the action on $M$ of a subgroup $G_{A_c'} \subset \mathcal{P}(M)$ that we will construct later on by finding a completion $\mathcal{A}_c'$ of the standard polar family $\mathcal{A}_G'$. This goal, that in principle seems rather technical, gains importance when we recall the definition of the symplectic leaves (Definition 4.1) and the Symplectic Leaf Correspondence Theorem (Theorem 4.5) where we saw that all these ideas are well behaved when we deal with the quotients of $M$ by *genuine subgroups* of $\mathcal{P}(M)$.

As a corollary we will obtain a correspondence between the symplectic leaves of the two legs of the diagram $(M/A_G, \{\cdot, \cdot\}_{M/A_G}) \xleftarrow{\pi} (M, \{\cdot, \cdot\}) \xrightarrow{\vartheta} (M/A'_G, \{\cdot, \cdot\}_{M/A'_G})$ in many of the situations identified in the preceding paragraphs in which that diagram is a dual pair.

**Proposition 7.10** *Let $G$ be a Lie group that acts canonically and properly on the symplectic manifold $(M, \omega)$ and $A_G \subset \mathcal{P}(M)$ be the associated Poisson diffeomorphisms subgroup. Then, the standard polar family of vector fields $\mathcal{A}'_G$ admits a completion $\mathcal{A}'_c$ that makes the polar distribution $A'_G$ complete.*

**Proof.**     The main tool in the proof will be the reconstruction equations presented in Appendix 8.2, hence the reader interested in the presentation that follows is encouraged to make a forward excursion to that section in order to get acquainted with the notation that we will use in the following paragraphs without much explanation.

Let $m \in M$ be an arbitrary point. Theorem 8.1 guarantees the existence of a $G$–invariant neighborhood $U$ of $m$ and of a $G$–equivariant symplectomorphism $\phi : U \to Y_r := G \times_{G_m} (\mathfrak{m}_r^* \times (V_m)_r)$ satisfying $\phi(m) = [e, 0, 0]$. Since the reconstruction equations (8.5) through (8.7) provide us with an explicit way to write down the Hamiltonian vector field associated to any function $h \in C^\infty(Y_r)^G$ we will use them to find a suitable generating family of complete vector fields for the polar distribution $A'_G$ by working in all the possible tubes $Y_r$ and translating our results back to $M$ via the symplectic diffeomorphisms $\phi$. The following arguments explain in detail this strategy.

Let $C_c^\infty(W)^{G_m}$, with $W \subset \mathfrak{m}_r^* \times (V_m)_r$ an open $G_m$–invariant neighborhood of the origin, be the set of compactly supported $G_m$–invariant smooth functions on $W$. Since the subgroup $G_m$ is compact and fixes the origin of $\mathfrak{m}_r^* \times (V_m)_r$, it is clear that

$$\{\mathbf{d}h(0, 0) \mid h \in C_c^\infty(W)^{G_m}\} = \{\mathbf{d}f(0, 0) \mid f \in C^\infty(\mathfrak{m}_r^* \times (V_m)_r)^{G_m}\},$$

and therefore it suffices to use the set of Hamiltonian vector fields associated to the functions in $C_c^\infty(W)^{G_m}$ to generate the image by $\phi$ of the polar distribution $A'_G(m)$ evaluated at $m$. By inspection of the reconstruction equations it is easy to see that any of those Hamiltonian vector fields is complete: take for instance (8.7), with $h \circ \pi \in C_c^\infty(W)^{G_m}$ and consider it as a vector field on $(V_m)_r$ with parameters $g$ and $\rho$. The compact support condition on $h \circ \pi$ implies that for any value of the parameters $g$ and $\rho$ the vector field $X_{V_m}(g, \rho, v)$ on $(V_m)_r$ is compactly supported and therefore complete. A similar analysis on (8.5) and (8.6) reveals the same conclusions for $X_G$ and $X_{\mathfrak{m}^*}$, which proves that $X_h$ is a complete vector field on $Y_r$. Also, since $\phi$ is a symplectic map, the vector field $X_{h \circ \phi} = T\phi^{-1} \circ X_h \circ \phi$ on $U$ is also complete and as it is zero outside $\phi^{-1}(G \times_{G_m} W) \subset U$ it can be trivially extended to a complete Hamiltonian vector field on $M$ associated to a $G$–invariant function (any extension of $h \circ \phi$). The union of all the similarly constructed vector fields using as many tubular neighborhoods as necessary constitutes a completion of the standard polar family $\mathcal{A}'_G$ hence proving that $A'_G$ is complete.    ■

The proposition that we just proved guarantees that the symplectic leaves of $M/A'_G$ are well defined for proper symplectic actions. Moreover, the combination of this statement with the Symplectic Leaf Correspondence Theorem 4.5, as well as with Theorems 7.2, 7.6, and 7.9, guarantees the correspondence between the symplectic leaves on the legs of the diagram $(M/A_G, \{\cdot, \cdot\}_{M/A_G}) \xleftarrow{\pi} (M, \omega) \xrightarrow{\vartheta} (M/A'_G, \{\cdot, \cdot\}_{M/A'_G})$ in a variety of situations that we enumerate for completeness in the following corollary.

**Corollary 7.11** *Let $G$ be a Lie group that acts canonically and properly on the symplectic manifold $(M, \omega)$ and $A_G \subset \mathcal{P}(M)$ be the associated Poisson diffeomorphisms subgroup. Suppose that at least ONE of the following conditions is satisfied:*

**(i)** *$A_G$ is (weakly or strongly) Hamiltonian,*

**(ii)** *the Lie group $G$ is compact and connected, and the tori of $A_G$ are coisotropic,*

**(iii)** *for any point $m \in M$, the orbit $G \cdot m$ is coisotropic, there exists a $\mathrm{Ad}(K_m)$-invariant complement to $\mathfrak{k}_m$ in $\mathfrak{g}$, $\mathfrak{k}_m$ is Abelian, and $\gamma_m$ is exact (see Theorem 7.9 for this notation).*

*Then, the map*

$$(M/A_G)/\{\cdot, \cdot\}_{M/A_G} \quad \longrightarrow \quad (M/A'_G)/\{\cdot, \cdot\}_{M/A'_G}$$
$$\mathcal{L}_{[m]_{A_G}} \quad \longmapsto \quad \mathcal{L}_{[m]_{A'_G}}$$

*establishes a bijection between the symplectic leaves of $(M/A_G, \{\cdot, \cdot\}_{M/A_G})$ and those of $(M/A'_G, \{\cdot, \cdot\}_{M/A'_G})$*

**Remark 7.12** As we recalled in Remark 7.12 the symplectic leaves of the Poisson manifolds in the legs of a Lie–Weinstein dual pair $(P_1, \{\cdot, \cdot\}_{P_1}) \xleftarrow{\pi_1} (M, \omega) \xrightarrow{\pi_2} (P_2, \{\cdot, \cdot\}_{P_2})$ in which the projections $\pi_1$ and $\pi_2$ have connected fibers are in bijection. Moreover, it can be shown [Bl01] that if $\mathcal{L}_1$ and $\mathcal{L}_2$ are symplectic leaves of $P_1$ and $P_2$ in correspondence and $\mathcal{K} \subset M$ is the immersed connected submanifold of $M$ such that $\mathcal{K} = \pi_1^{-1}(\mathcal{L}_1) = \pi_2^{-1}(\mathcal{L}_2)$ then, the symplectic forms $\omega_{\mathcal{L}_1}$ and $\omega_{\mathcal{L}_2}$ on $\mathcal{L}_1$ and $\mathcal{L}_2$, respectively, are related by the equation

$$i_{\mathcal{K}}^* \omega = \pi_1|_{\mathcal{K}}^* \omega_{\mathcal{L}_1} + \pi_2|_{\mathcal{K}}^* \omega_{\mathcal{L}_2}, \tag{7.2}$$

where $i_{\mathcal{K}} : \mathcal{K} \hookrightarrow M$ denotes the natural inclusion. In my forthcoming paper [O02] I will show that the symplectic leaves in the von Neumann pairs studied in this section can be, under certain hypotheses, be endowed with actual smooth structures that make them into real symplectic manifolds related by an equality identical to (7.2).   ♦

# 8   Appendices

In this section we explain more in detail some of the tools that have been used throughout the paper.

## 8.1   A normal form for canonical proper actions

The Slice Theorem in the category of globally Hamiltonian proper Lie group actions (in this section the expression globally Hamiltonian means that the action has a globally defined equivariant momentum map associated) is a well known and widely used tool introduced by Marle [Mar84, Mar85] and by Guillemin and Sternberg [GS84a]. The classical construction of the normal form coordinates in this setup uses very strongly the existence of a momentum map in the space that we want to locally model. In the following paragraphs we show how to reproduce this result for canonical proper actions on a symplectic manifold that do not necessarily have a momentum map associated. In the exposition we will limit ourselves to present the ingredients of this construction. For a complete presentation the reader is encouraged to check with [OR01] where this normal form is explained in full detail.

Let $(M, \omega)$ be a symplectic manifold, $G$ be a Lie group acting properly and canonically on it, and $m \in M$ be an arbitrary point in $M$ around which we want to construct the slice coordinates. The following facts can be readily verified:

**(i)** The vector space $V_m := T_m(G \cdot m)^\omega/(T_m(G \cdot m)^\omega \cap T_m(G \cdot m))$ is symplectic with the symplectic form $\omega_{V_m}$ defined by

$$\omega_{V_m}([v], [w]) := \omega(m)(v, w),$$

for any $[v] = \pi(v)$ and $[w] = \pi(w) \in V_m$, and where $\pi : T_m(G \cdot m)^\omega \to T_m(G \cdot m)^\omega/(T_m(G \cdot m)^\omega \cap T_m(G \cdot m))$ is the canonical projection. The vector space $V_m$ is called the **symplectic normal space** at $m$.

**(ii)** Let $H := G_m$ be the isotropy subgroup of $m$. The properness of the $G$–action guarantees that $H$ is compact. The mapping $(h, [v]) \longmapsto [h \cdot v]$, with $h \in H$ and $[v] \in V_m$, defines a linear canonical action of the Lie group $H$ on $(V_m, \omega_{V_m})$, where $g \cdot u$ denotes the tangent lift of the $G$–action on $TM$, for $g \in G$ and $u \in TM$. We will denote by $\mathbf{J}_{V_m} : V_m \longrightarrow \mathfrak{h}^*$ the associated $H$-equivariant momentum map.

**(iii)** Let $\mathfrak{g}$ and $\mathfrak{h}$ be the Lie algebras of $G$ and $H$, respectively. The vector subspace $\mathfrak{k} \subset \mathfrak{g}$ given by

$$\mathfrak{k} = \{\eta \in \mathfrak{g} \mid \eta_M(m) \in (\mathfrak{g} \cdot m)^\omega\}$$

is a subalgebra of $\mathfrak{g}$ such that $\mathfrak{h} \subset \mathfrak{k}$.

**(iv)** Lie algebra decompositions: the compactness of the isotropy subgroup $H$ allows us to choose two $\mathrm{Ad}_H$-invariant complements: $\mathfrak{m}$ to $\mathfrak{h}$ in $\mathfrak{k}$ and $\mathfrak{q}$ to $\mathfrak{k}$ in $\mathfrak{g}$. Therefore, we have the orthogonal decompositions

$$\mathfrak{g} = \mathfrak{h} \oplus \mathfrak{m} \oplus \mathfrak{q}, \quad \text{where} \quad \mathfrak{k} = \mathfrak{h} \oplus \mathfrak{m}, \tag{8.1}$$

as well as their duals

$$\mathfrak{g}^* = \mathfrak{h}^* \oplus \mathfrak{m}^* \oplus \mathfrak{q}^*, \quad \text{where} \quad \mathfrak{k}^* = \mathfrak{h}^* \oplus \mathfrak{m}^*.$$

**(v)** **The normal form tube**: there are two $H$-invariant neighborhoods $(V_m)_r$ and $\mathfrak{m}_r^*$ of the origin in $V_m$ and $\mathfrak{m}^*$ such that the twisted product

$$Y_r := G \times_H (\mathfrak{m}_r^* \times (V_m)_r) \tag{8.2}$$

is a symplectic manifold acted on by $G$ according to the expression $g \cdot [h, \eta, v] := [gh, \eta, v]$, for any $g \in G$, and any $[h, \eta, v] \in Y_r$. This action is canonical. The use of the construction (8.2) is justified by the following theorem.

**Theorem 8.1 (Slice Theorem for canonical Lie group actions)** *Let $(M, \omega)$ be a symplectic manifold and let $G$ be a Lie group acting properly and canonically on $M$. Let $m \in M$ and denote $H := G_m$. Then, the manifold*

$$Y_r := G \times_H (\mathfrak{m}_r^* \times (V_m)_r) \tag{8.3}$$

*introduced in (8.2) is a symplectic $G$–space and can be chosen such that there is a $G$–invariant neighborhood $U$ of $m$ in $M$ and an equivariant symplectomorphism $\phi : U \to Y_r$ satisfying $\phi(m) = [e, 0, 0]$.*

## 8.2 The reconstruction equations for canonical proper Lie group actions

The reconstruction equations are the differential equations that determine the Hamiltonian vector field associated to a $G$–invariant Hamiltonian in the coordinates provided by Theorem 8.1. In the globally Hamiltonian context these equations can be found in [O98, RWL99, OR02a]. As we will see, it is remarkable that in the absence of a momentum map, the reconstruction equations written using the symplectic form of (8.2) are *formally* identical to the ones obtained for the globally Hamiltonian case.

In order to present the reconstruction equations let $h \in C^\infty(Y_r)^G$ be the Hamiltonian function whose associated vector field $X_h$ we want to write down. Let $\pi : G \times \mathfrak{m}_r^* \times (V_m)_r \to G \times_H (\mathfrak{m}_r^* \times (V_m)_r)$ be the canonical projection. The $G$-invariance of $h$ implies that $h \circ \pi \in C^\infty(G \times \mathfrak{m}_r^* \times (V_m)_r)^H$ can be understood as a $H$–invariant function that depends only on the $\mathfrak{m}_r^*$ and $(V_m)_r$ variables, that is,

$$h \circ \pi \in C^\infty(\mathfrak{m}_r^* \times (V_m)_r)^H.$$

Now, since the projection $\pi$ is a surjective submersion, the Hamiltonian vector field $X_h$ can be locally expressed as

$$X_h = T\pi(X_G, \, X_{\mathfrak{m}^*}, \, X_{V_m}),$$

with $X_G$, $X_{\mathfrak{m}^*}$ and $X_{V_m}$ locally defined smooth maps on $Y_r$ and having values in $TG$, $T\mathfrak{m}_r^*$ and $T(V_m)_r$ respectively. Moreover, using the $\mathrm{Ad}_H$–invariant decomposition of the Lie algebra $\mathfrak{g}$ introduced in the previous section, the mapping $X_G$ can be written, for any $(g, \rho, v) \in G \times \mathfrak{m}_r^* \times (V_m)_r$, as

$$X_G(g, \rho, v) = T_e L_g \big( X_\mathfrak{h}(g, \rho, v) + X_\mathfrak{m}(g, \rho, v) + X_\mathfrak{q}(g, \rho, v) \big),$$

with $X_\mathfrak{h}$, $X_\mathfrak{m}$, and $X_\mathfrak{q}$, locally defined smooth maps on $Y_r$ with values in $\mathfrak{h}$, $\mathfrak{m}$, and $\mathfrak{q}$, respectively. In what follows we give the expressions that determine $X_G$, $X_{\mathfrak{m}^*}$, and $X_{V_m}$ as a function of the differential of the Hamiltonian $h$.

First, the construction of $\mathfrak{q}$ as a complement to $\mathfrak{k}$ guarantees that the bilinear pairing $\langle \cdot, \cdot \rangle_\mathfrak{q}$ in $\mathfrak{q}$ defined by

$$\langle \xi, \eta \rangle_\mathfrak{q} := \omega(m) \left( \xi_M(m), \eta_M(m) \right)$$

is non degenerate. Let $\mathbb{P}_{\mathfrak{h}^*}$, $\mathbb{P}_{\mathfrak{m}^*}$, and $\mathbb{P}_{\mathfrak{q}^*}$ denote the projections from $\mathfrak{g}^*$ onto $\mathfrak{h}^*$, $\mathfrak{m}^*$, and $\mathfrak{q}^*$, respectively, according to the $\mathrm{Ad}_H^*$–invariant splitting $\mathfrak{g}^* = \mathfrak{h}^* \oplus \mathfrak{m}^* \oplus \mathfrak{q}^*$. The non degeneracy of $\langle \cdot, \cdot \rangle_\mathfrak{q}$ implies that the mapping

$$
\begin{array}{rccc}
F: & \mathfrak{k} \times \mathfrak{k}^* \times \mathfrak{q} & \longrightarrow & \mathfrak{q}^* \\
 & (\xi, \lambda, \tau) & \longmapsto & \mathbb{P}_{\mathfrak{q}^*} \left( \mathrm{ad}^*_{(\xi+\tau)} \lambda \right) + \langle \tau, \cdot \rangle_\mathfrak{q},
\end{array}
\tag{8.4}
$$

is such that $F(0,0,0) = 0$ and that its derivative $D_\mathfrak{q} F(0,0,0) : \mathfrak{q}^* \to \mathfrak{q}^*$ is a linear isomorphism. The Implicit Function Theorem implies the existence of a locally defined function $\tau : \mathfrak{k} \times \mathfrak{k}^* \to \mathfrak{q}$ around the origin such that $\tau(0,0) = 0$, and $\mathbb{P}_{\mathfrak{q}^*} \left( \mathrm{ad}^*_{(\xi+\tau(\xi,\lambda))} \lambda \right) + \langle \tau(\xi, \lambda), \cdot \rangle_\mathfrak{q} = 0$. Let now $\psi : \mathfrak{m}^* \times V_m \to \mathfrak{q}$ be the locally defined function given by:

$$\psi(\rho, v) := \tau \left( D_{\mathfrak{m}^*}(h \circ \pi)(\rho, v), \rho + \mathbf{J}_{V_m}(v) \right).$$

With these expressions at hand it can be verified [OR01] that $X_h$ is given by

$$X_G(g, \rho, v) = T_e L_g \left( \psi(\rho, v) + D_{\mathfrak{m}^*}(h \circ \pi)(\rho, v) \right) \tag{8.5}$$

$$X_{\mathfrak{m}^*}(g, \rho, v) = \mathbb{P}_{\mathfrak{m}^*} \left( \mathrm{ad}^*_{D_{\mathfrak{m}^*}(h \circ \pi)} \rho \right) + \mathrm{ad}^*_{D_{\mathfrak{m}^*}(h \circ \pi)} \mathbf{J}_{V_m}(v) + \mathbb{P}_{\mathfrak{m}^*} \left( \mathrm{ad}^*_{\psi(\rho, v)}(\rho + \mathbf{J}_{V_m}(v)) \right) \tag{8.6}$$

$$X_{V_m}(g, \rho, v) = B_{V_m}^\sharp \left( D_{V_m}(h \circ \pi)(\rho, v) \right) \tag{8.7}$$

## 8.3    Tubewise Hamiltonian actions

In the following paragraphs we will see how the normal form introduced in Section 8.1 helps us determining, in the framework of proper canonical actions, when such an action is tubewise Hamiltonian. In order to be more specific, recall that by Theorem 8.1, any $G$-orbit of a symplectic $G$-space, $(M, \omega)$ has an invariant neighborhood around the orbit $G \cdot m$ that can be modeled by an associated bundle like the one presented in (8.3). Consequently, we can conclude that the canonical proper $G$-action on $(M, \omega)$ is strongly (resp. weakly) tubewise Hamiltonian if the $G$-action on each $G$-invariant coordinate patch (8.3) is strongly (resp. weakly) Hamiltonian. The following result provides a sufficient condition regarding the strong case whose proof can be found in [OR01].

**Proposition 8.2** *Let $(M, \omega)$ be a symplectic manifold and let $G$ be a Lie group with Lie algebra $\mathfrak{g}$, acting properly and canonically on $M$. Let $m \in M$, $H := G_m$ and $Y_r := G \times_H (\mathfrak{m}_r^* \times (V_m)_r)$ be the normal coordinates around the orbit $G \cdot m$ introduced in (8.3). If the $G$-equivariant, $\mathfrak{g}^*$–valued one form $\gamma \in \Omega^1(G; \mathfrak{g}^*)$ defined by*

$$\langle \gamma(g) \cdot T_e L_g \cdot \eta, \xi \rangle := -\omega(m) \left( \left( \mathrm{Ad}_{g^{-1}} \xi \right)_M (m), \eta_M(m) \right) \quad \text{for any} \quad g \in G, \, \xi, \eta \in \mathfrak{g} \tag{8.8}$$

*is exact then, the $G$–action on $Y_r$ given by $g \cdot [h, \eta, v] := [gh, \eta, v]$, for any $g \in G$, and any $[h, \eta, v] \in Y_r$, has a standard momentum map associated. Therefore, if the group $G$ is connected and (8.8) is exact, the $G$–action on $Y_r$ is strongly Hamiltonian.*

The next proposition provides another characterization of the exactness of (8.8) and therefore another sufficient condition for the $G$–action on the tube $Y_r$ to be strongly Hamiltonian when $G$ is connected. See [OR01] for a proof.

**Proposition 8.3** *Suppose that we are in the conditions of Proposition 8.2. Let $m \in M$, $H := G_m$ and $Y_r := G \times_H (\mathfrak{m}_r^* \times (V_m)_r)$ be the normal coordinates around the orbit $G \cdot m$ introduced in (8.3). let $\Sigma : \mathfrak{g} \times \mathfrak{g} \longrightarrow \mathbb{R}$ be the two cocycle given by*

$$\Sigma(\xi, \eta) = \omega(m) \left( \xi_M(m), \eta_M(m) \right), \quad \xi, \eta \in \mathfrak{g},$$

*and $\Sigma^\flat : \mathfrak{g} \to \mathfrak{g}^*$ be defined as $\Sigma^\flat(\xi) = \Sigma(\xi, \cdot)$, $\xi \in \mathfrak{g}$. Then, the form (8.8) is exact if and only if there exists a $\mathfrak{g}^*$–valued group one cocycle $\theta : G \to \mathfrak{g}^*$ such that*

$$T_e \theta = \Sigma^\flat.$$

*In the presence of this cocycle, the map $\mathbf{J}_\theta : G \times \mathfrak{k}^* \to \mathfrak{g}^*$ given by*

$$\mathbf{J}_\theta(g, \nu) := \mathrm{Ad}^*_{g^{-1}} \nu - \theta(g) \tag{8.9}$$

*is a momentum map for the $G$–action on $G \times \mathfrak{k}^*$ with non equivariance cocycle equal to $-\theta$.*

The following corollary presents two situations in which the hypotheses of Proposition 8.2 are trivially satisfied for any point $m \in M$.

**Corollary 8.4** *Let $(M, \omega)$ be a symplectic manifold and let $G$ be a connected Lie group with Lie algebra $\mathfrak{g}$, acting properly and canonically on $M$. If either,*

**(i)** $H^1(G) = 0$, *or*

**(ii)** *all the $G$–orbits are isotropic*

*then, the associated subgroup $A_G$ of $\mathcal{P}(M)$ is tubewise strongly Hamiltonian.*

## 8.4 Proof of Lemma 3.6

Let $D$ be a smooth integrable regular distribution on the manifold $M$ and $\pi_D : M \to M/D$ be the associated surjective submersion. Let $U$ be a $D$–saturated open subset of $M$, $z$ a point in $U$, and $f \in C^\infty(U)^D$. Since $\pi_D$ is a submersion, there are charts $(V, \varphi)$ and $(W, \psi)$ around $z$ and $\pi_D(z)$, respectively, such that $\pi_D(V) = W$, $\varphi : V \to V' \times W'$, $\psi : W \to W'$, $\varphi(z) = (0, 0)$, and the local representative of $\pi_D$, that is, $\psi \circ \pi_D \circ \varphi^{-1} : V' \times W' \to W'$ is the projection onto the second factor. We will shrink $V$ if necessary so that $V \subset U$.

Let $B_\epsilon(0) \subset W'$ be a ball of radius $\epsilon > 0$ and $\phi : B_\epsilon(0) \to [0, 1]$ be a bump function such that:

$$\phi|_{B_{\epsilon/2}(0)} = 1 \quad \text{and} \quad \phi|_{B_\epsilon(0) \backslash B_{3\epsilon/4}(0)} = 0.$$

Let $F'_D : \psi^{-1}(B_\epsilon(0)) \subset W \to \mathbb{R}$ be the smooth function given by $F'_D(l) = f_D(l)\phi(\psi(l))$, $l \in \psi^{-1}(B_\epsilon(0))$, where $f_D : \pi_D(U) \to \mathbb{R}$ is the unique smooth function determined by the relation $f = f_D \circ \pi_D|_U$. As the function $F'_D$ and all its derivatives are zero in the boundary of $\psi^{-1}(B_\epsilon(0))$, $F'_D$ can be extended to a function $F_D \in C^\infty(M/D)$.

Let $F \in C^\infty(M)^D$ be the function given by $F := F_D \circ \pi_D$ and $\Sigma$ be the submanifold of $M$ through $z$ defined as $\Sigma := \varphi^{-1}(\{0\} \times B_{\epsilon/2}(0))$. Notice that $\pi_D(\Sigma)$ is an open subset of $\pi_D(U)$ since $\psi(\pi_D(\Sigma)) = B_{\epsilon/2}(0)$ is an open subset of $\psi(\pi_D(U))$. Consequently, $T = \pi_D^{-1}(\pi_D(\Sigma)) \subset U$ is an open $D$–invariant subset of $U$.

We will prove the lemma by showing that $F|_T = f|_T$. Indeed, let $m \in T$ arbitrary. By definition $m \in T$ iff $\pi_D(m) \in \pi_D(\Sigma)$ or, equivalently, there exists an element $z \in \varphi^{-1}(\{0\} \times B_{\epsilon/2}(0))$ such that $\pi_D(m) = \pi_D(z)$. Due to the local expression of $\pi_D$ in the charts $(V, \varphi)$ and $(W, \psi)$ we have that $\psi(\pi_D(z)) \in B_{\epsilon/2}(0)$ or, equivalently, $\pi_D(z) \in \psi^{-1}(B_{\epsilon/2}(0))$. With this in mind, we have that

$$F(m) = F_D \circ \pi_D(m) = F_D \circ \pi_D(z) = F'_D \circ \pi_D(z),$$

where the previous equality follows from the fact that $\pi_D(z) \in \psi^{-1}(B_{\epsilon/2}(0))$. We now use the definition of $F'_D$ and

$$F'_D \circ \pi_D(z) = f_D(\pi_D(z))\phi(\psi(\pi_D(z))) = f_D(\pi_D(z)) = f_D(\pi_D(m)) = f(m),$$

which proves that $F(m) = f(m)$, as required. ∎

**Acknowledgments** I thank Rui Loja Fernandes for his interest in this project and for his numerous and relevant comments to the first draft. I am grateful for useful and interesting conversations with James Montaldi and Tudor Ratiu. I also thank the hospitality of Jerry Marsden during my stay in Caltech in March 2001, where the main ideas of this paper were gestated. This research was partially supported by the European Commission through funding for the Research Training Network *Mechanics and Symmetry in Europe* (MASIE).